\documentclass[12pt]{article}
\usepackage{amssymb,amsthm,caption2,vmargin}

\input epsf

\setpapersize{USletter}

\setmarginsrb{1.4in}{1.5in}{1.4in}{1in} {0pt}{1mm}{1em}{3em}
\usepackage{setspace}

\usepackage{amsmath,amsfonts}
\usepackage[mathscr]{eucal}
\usepackage{amssymb,amsmath}
\usepackage{latexsym}
\usepackage{amsfonts}
\usepackage{amscd}
\usepackage{graphicx}
\usepackage{mathptmx}
\newcommand {\bt}{\tilde{\beta}}

\newcommand {\R}{\mathbb {R}}
\newcommand {\T}{\mathbb {T}}
\newcommand {\cH}{\mathcal {H}}
\newcommand {\G}{\mathcal {G}}

\newcommand {\D}{\mathcal {D}}

\newcommand {\F}{\mathcal {F}}
\newcommand {\N}{\mathcal {N}}
\newcommand {\Z}{\mathbb {Z}}
\newcommand {\half}{\frac {1}{2}}
\newcommand {\third}{\frac {1}{3}}
\newcommand {\quarter}{\frac {1}{4}}

\newcommand {\eighth}{\frac {1}{8}}
\newcommand {\tenth}{\frac {1}{10}}
\newcommand {\cfour}{{xx}}
\newcommand {\ceight}{{x}}
\newcommand {\csixtyfour}{{}}
\newcommand {\st}{\text{ such that }}

\newcommand {\then}{\Rightarrow}

\newcommand {\suml}{\sum\limits}

\newcommand {\supl}{\sup\limits}

\newcommand{\norm}[1]{\left\lVert{#1}\right\rVert}
\newcommand{\abs}[1]{\lvert{#1}\rvert}
\newcommand{\diam}{{\rm diam}}
\newcommand{\radius}{{\rm radius}}

\newcommand{\dist}{{\rm dist}}
\newcommand{\ball}{{\rm Ball}}

\newcommand{\supp}{{\rm supp}}
\newcommand{\length}{\ell}
\newtheorem{theorem}{Theorem}[section]
\newtheorem{lemma}[theorem]{Lemma}
\newtheorem{proposition}[theorem]{Proposition}
\newtheorem{cor}[theorem]{Corollary}
\theoremstyle{remark}
\newtheorem{rem}[theorem]{Remark}

\numberwithin{equation}{section}

\begin{document}

\title{Subsets of Rectifiable curves in Hilbert Space-The Analyst's TSP}
\author{Raanan Schul\\
\tt{schul@math.ucla.edu}}         
\date{}
\maketitle\thispagestyle{empty}

\begin{abstract}

We study one dimensional sets (Hausdorff dimension) lying in a Hilbert space.  The aim is to classify subsets of  Hilbert spaces that are contained in a connected set of finite Hausdorff length. We do so by extending  and improving results of Peter Jones and Kate Okikiolu for sets in $\R^d$. Their results formed the basis of quantitative rectifiability in $\R^d$.   
We prove  a quantitative version of the following statement:  a connected set of finite Hausdorff  length  (or a subset of one), is characterized by the fact that inside balls at most scales around most points of the set, the set lies close to a straight line segment (which depends on the ball).
This is done via a quantity, similar to the one introduced in \cite{J1}, which is a geometric analog of the  Square function.  This allows us to  conclude that for a given set $K$,  the $\ell_2$ norm of this  quantity (which is a function of  $K$) has size comparable to a shortest (Hausdorff length) connected set containing $K$. 
In particular, our results imply that, with a correct reformulation of the theorems, the estimates in \cite{J1,Ok} are {\bf  independent of the ambient dimension}.

\noindent
{\bf Mathematics Subject Classification (2000):} 28A75
\end{abstract}

\tableofcontents

\section{Introduction}
\label{intro}
%
\subsection{Basic Notation and Definitions}
We start with some basic definitions and some  history.  We will state our new results in section \ref{sect_new_results}.
\subsubsection*{Cubes, Grids, Balls, and Nets.  Multiresolution Families}
A cube $Q$ in $\R^d$ is a set of the form $I_1\times I_2 \times ... \times I_d$, where
$I_1,I_2,...,I_d$ are intervals satisfying $\abs{I_1}=\abs{I_2}=...=\abs{I_d}=l(Q)$ .
We call $l(Q)$  the side-length of $Q$.
We denote by $\lambda Q$ the cube with the same center as $Q$, but with side-length   $\lambda l(Q)$.

A dyadic cube is a cube of the form
\begin{gather*}
Q=[{i_1\over 2^j},{i_1 +1\over 2^j}]\times...\times[{i_d\over 2^j},{i_d +1\over 2^j}]
\end{gather*}
where  $i_1,...,i_d,j$ are integers.
The standard dyadic grid on $\R^d$ is 
\begin{gather*}
	\D=\{Q=[{i_1\over 2^j},{i_1 +1\over 2^j}]\times...\times[{i_d\over 2^j},{i_d +1\over 2^j}]:
         		               i_1,...,i_d,j \text{ integers}\}.
\end{gather*}

A ball $Q$ is a set 
\begin{gather*}
\ball(x,r):=\{y:\norm{y-x}\leq r\}.
\end{gather*}
We denote by $l(Q)$ the diameter of the ball $Q$
and by $\lambda Q$ the ball with the same center as Q, but with radius $\lambda r$ instead of $r$ (we call this a dilation by $\lambda$ of $Q$).

We say that $X$ is an $\epsilon-net$ for $K$ if
\begin{quote} 
(i) $X\subset K$\\
(ii) $\norm{x_1-x_2} > \epsilon,\forall x_1,x_2\in X$\\
(iii)$\forall y\in K, \exists x\in X$ such that $\norm{x-y}\leq \epsilon$
\end{quote}
Hence $K\subset \bigcup\limits_{x\in X}\ball(x,\epsilon)$.  
Note that if $X' \subset K$ satisfies $\norm{x_1-x_2} > \epsilon,\forall x_1,x_2\in X'$ then $X'$ can be extended to an  $\epsilon-net$ $X$ since a maximal subset of $K$ satisfying (ii), will satisfy (iii).

Fix a set $K$. Denote by $X^K_n$ a  sequence of $2^{-n}-nets$ for $K$, such that $X^K_n \subset X^K_{n+1}$.
Set
\begin{gather}\label{G-K-def}
\hat{\G}^K=\{\ball(x,A2^{-n}):x\in X^K_n, n \text{ an integer}, n\geq n_0\}
\end{gather}
for a constant $A>1$ and $n_0$ an arbitrary (possibly negative) integer.  Existence of such a sequence of nets is assured since we may start by choosing a maximal subset of $K$ satisfying (ii) for  $n_0$ and then proceed inductively for $n>n_0$. (This is the only use of $n_0$. Unless explicitly stated, all results will be independent of $n_0$ and hence we will suppress it in the notation.)   

We call $\hat{\G}^K$ a {\it multiresolution family}.
Note that $\hat{\G}^K$ depends on $K$.
We also call  the standard dyadic grid a {\it multiresolution family}.

For a multiresolution $\hat{\G}$, we denote by $\lambda\hat{\G}$ the multiresolution given by dilating each element in $\hat{\G}$ by $\lambda$.
\subsubsection*{Neighborhoods}
We denote the $\epsilon$ neighborhood of a set $E$ by $\N_\epsilon(E)$.
\subsubsection*{Hausdorff Length and Arclength}
For a set $K$ we denote by  $\cH^1(K)$ the one dimensional  Hausdorff measure, which we call  \textit{Hausdorff length}.   See \cite{Ma} for definition and discussion.  
For a Lipschitz function (see below)  $\tau:[a,b]\to H$ (a Hilbert space) we will  denote by $\length(\tau)$ the arclength of $\tau$. We will also extend this definition to Borel sets of the domain of a given Lipschitz function and use it for the push-forward of this measure.    
%
\subsubsection*{Hilbert Space}
We shall concern ourselves with subsets of a Hilbert space which are subsets of finite length connected sets.  
All finite  length sets  are separable. 
Hence we shall only concern ourselves with separable subspaces of Hilbert spaces, which in turn, are separable Hilbert spaces.  Those are all isometric to subspaces of $\ell_2$ as vector spaces.  Note that since the isometries in question are of vector spaces, straight lines go to straight lines.  This will be  crucial so that we do not loose generality.  Hence we restrict our discussion to separable sets $K$ and  fix  $H=\ell_2$ as our Hilbert space.
\subsubsection*{Lipschitz Functions, Rectifiable Sets, Rectifiable curves}
A function $f:\R^k \to H$  is said to be \textit{Lipschitz} if
\begin{gather*}
{\norm{f(x)-f(y)} \over \norm{x-y}} \leq C_f,\forall x,y\in \R^k.
\end{gather*}
A set is called \textit{k-rectifiable} if it is contained in a countable union of images 
of  Lipschitz functions $f_j:\R^k \to H$, except for a set of k-dimensional Hausdorff measure zero.
For more details see \cite{Ma}, where one can also find an excellent discussion of rectifiability in the setting of $\R^d$, part of which carries over to the setting of $H$. 

A set is called a \textit{rectifiable curve} if it is the image a Lipschitz function defined on $\R$.
\subsubsection*{The Jones $\beta$ Numbers.  The Jones Function}
Assume we have a set $K$ lying in $\R^d$ or $H$.
Consider $Q$  a cube or ball.
We define the Jones $\beta$ number as 
\begin{eqnarray*}
\beta_{K}(Q)&=&\frac{2}{\diam(Q)} \inf\limits_{L \text{ line}}\sup\limits_{x \in K \cap Q }\dist(x,L)\\
&=&{\text{width of thinnest cylinder containing }K \cap Q \over \diam(Q)}.
\end{eqnarray*}
\begin{figure}[h]
\begin{center}
\scalebox{0.4}{\includegraphics*[-0in,-1in][8in,4in]{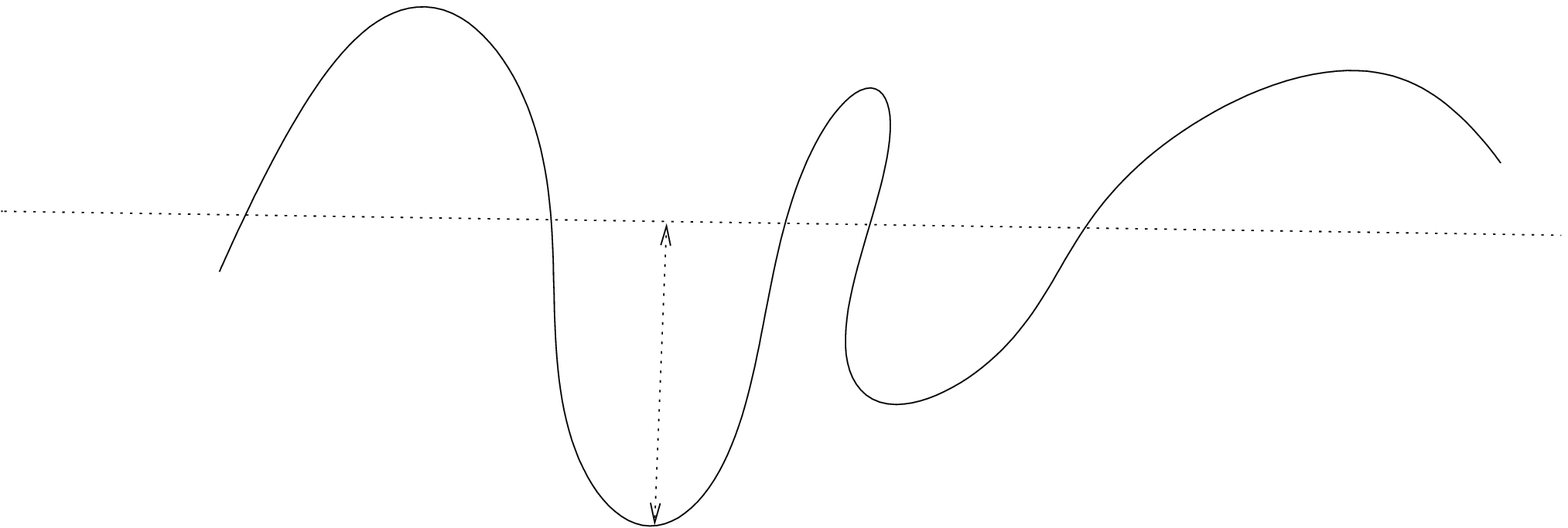}}
\put(-140,50){\tiny $D$}
\put(-60,100){\tiny $K\cap Q$}
\put(-220,80){\tiny $L$}
\end{center}
\caption[$\beta_K(Q)$]{$D=\half\beta_{K}(Q)diam(Q)$}
\label{beta}
\end{figure}

Hence if $K' \supset K$ then $\beta_{K'}(Q) \geq \beta_K(Q)$.
Note that we have defined a quantity which is scale independent.
This quantity is usually referred to as the Jones $\beta_\infty$ number in order to differentiate it from its $L^p$ variants (extensively developed by David and Semmes in \cite{DS} and generalized in \cite{Le}).  We omit the $\infty$ subscript as we will always use $\beta_\infty$.
We will occasionally use the notation $\beta_K(x,r):=\beta_K(\ball(x,r))$.
We will often omit $K$ from the notation when it is obvious what it is.

Fix a multiresolution family.
We define the Jones function $J(x)$ as follows:
\begin{gather*}
J(x)=\suml\beta^2(Q)\chi_Q(x)
\end{gather*}
where we sum over the  multiresolution family of $Q$'s we have fixed, and $\chi_Q$ is the indicator function of $Q$.  This should be thought of as an analog of (the square of)  the Square function for certain categories of sets.  In many cases one adds for each $Q$ a weight in the above sum. See \cite{BJ,DS,J1,Le,Ok} for explicit and implicit appearances.  One can view the remainder of this essay as an explanation of the right way of generalizing this notion to a larger category of sets. See subsection \ref{sect_new_results} for more details and more precise statements. 

\subsection{Subsets of Rectifiable Curves.  The Analyst's Traveling Salesman Problem}
\subsubsection*{Overview}
Given a set $K \subset H$ one can ask under what conditions is $K$ contained in $\Gamma$, 
the image of a single
Lipschitz function $\gamma:[0,1] \to H$. 
One can also ask for estimates on the minimal arc-length of such $\gamma$.    
(Recall that, up to multiplicative constants, the arclength of a Lipschitz curve $\gamma$ is equivalent to the 1-dimensional Hausdorff length of the image of $\gamma$. See Lemma \ref{parametrization})
In \cite{J1} Peter Jones gave such an estimate for a {\it planar} set $K$ in the form of an $\ell_2$ sum. 
In fact, he gave  a necessary and sufficient 
condition for a {\it planar} set $K$ to be contained in the image of a Lipschitz function by  requiring this $\ell_2$ sum to be finite.  He also gave a  construction of  
a curve whose image $\Gamma_0$ contains $K$, such that 
$\cH^1(\Gamma_0) \lesssim \cH^1 (\Gamma_{MST})$.
Here $\lesssim$ means `less then a constant multiple of', $\cH^1(\cdot)$ is the one dimensional Hausdorff measure  and   $\Gamma_{MST}$  is a shortest connected set containing $K$, whose existence is assured using Arzela-Ascoli and Golab's Theorem.   
The curve he constructed enjoys some useful  properties (for example, see \cite{BJ}).
Philosophically, one should view this $\ell_2$ sum as the expectation of  how much $K$ deviates from being flat (along a line segment) in a random window. The surprise was that this quantity ended up being equivalent to the length of the shortest curve containing $K$.
In \cite{Ok} Okikiolu extended this result to sets $K\subset \R^d$, rather then $\R^2$.  These results formed the basis of a theory now called `quantitative rectifiability' which was extended by many authors, of which we should give special mention to Guy David and Stephen Semmes whose work on Uniform Rectifiability inspired part of this essay. 
For example  see \cite{DS,BJ,Le}.  Also see \cite{Pa1} for a more complete survey and bibliography.  
As it turns out, many aspects of the quantitative form in which things will be presented are \textit{parallel} to the theory of wavelets, and a dictionary (discovered by Peter Jones) can be written (see Appendix B in \cite{Schul}). 

Unfortunately, the dependence of the constants in \cite{J1,Ok} on $d$ (as in $\R^d$) is exponential. 
This gives motivation to have a  Hilbert space version of this theory, which is equivalent to obtaining dimension free estimates.  This is the goal of this essay. 
One should note other examples in harmonic analysis where one was able to obtain dimension free estimates such as  the boundedness of the ball Maximal function and the norm of the size of the Riesz vector.  See \cite{St,StSt}.

A very natural question to ask is `How does this relate to the Euclidean TSP or Euclidean MST?' 
(the classical TSP  is finding a shortest Hamiltonian cycle on a finite graph;
the classical MST problem is finding a minimal spanning tree in a finite graph;  
their Euclidean counterparts are when the graphs are embedded in Euclidean space).
If one wants a polynomial time algorithm and is willing to accept an answer that is not `the shortest',
but `the shortest up to a constant multiple' then there are readily available algorithms (see e.g. \cite{JM} for a description of several algorithms).  They either do not come with a multiresolution analysis (such as a greedy algorithm which gives multiplicative constant $2$), or have exponential (super-exponential!) dependence on the dimension $d$ (such as \cite{Ar} which gives multiplicative constant $1+ \epsilon$).  
The results of \cite{J1}, \cite{Ok} can  be used (and are used) to give such algorithms, but with exponentially bad dependence of the constant on the dimension of the ambient space.  As an example of an application of our theorem  we give a proof that a {\it local  version} of the Farthest Insertion algortithm for  the MST converges to a connected set no longer than a constant multiple of the length of the MST. 
Our proof gives constants independent of the ambient dimension $d$!
One should note that the constants given by our proof are not as good as the ones experimentally found, as discussed in \cite{JM}.\\

The results of \cite{J1,Ok} can also be used to prove results regarding existence of Spanning Trees for rectifiable curves.  In particular they can be used to give an alternative geometric construction to  \cite{KK} that works in $\R^d$.  We have not yet been able to extend this result to the setting of a Hilbert space.

\subsubsection*{More Details}
We now discuss the results of \cite{J1,Ok} in a little more detail.
 
Jones (\cite{J1})  proved that for any curve $\gamma$ with image $\Gamma\subset \R^2$ (or equivalently, for any connected set $\Gamma\subset \R^2$)
\begin{gather}\label{J_sum_less_length}
	\suml_\D \beta_\Gamma^2(3Q)l(Q) \lesssim
	\cH^1(\Gamma).
\end{gather}
where $\D$ is the dyadic grid on $\R^2$ and $l(Q)$ is the side length of a cube $Q\in \D$.

Jones also gave a construction
that, given a set $K\subset\R^d$ (and in particular $\R^2$), yields a connected set 
$\Gamma_0\supset K$.     
The length of $\Gamma_0$ satisfies 
\begin{gather*}
	\cH^1(\Gamma_0) \lesssim 
	\diam(K) + \suml_\D \beta_K^2(3Q)l(Q)
\end{gather*}
where $\D$
is a dyadic grid on $\R^d$.
This construction is a multi-scale
algorithm, starting from the "roughest" scale and then refining.
This multi-scale method
also allows one to form approximations of the final connected set by applying only
a finite number of iterations.  
Combining this length estimate with  \eqref{J_sum_less_length} one gets
\begin{gather*}
\cH^1(\Gamma_0) \lesssim \cH^1(\Gamma_{MST})
\end{gather*}   
and
\begin{gather*}
\diam(K) + \suml_\D \beta_K^2(3Q)l(Q) \sim \cH^1(\Gamma_{MST})
\end{gather*}
as $\beta_K \leq \beta_{\Gamma_{MST}}$ and $\diam(K) \leq \cH^1(\Gamma_{MST})$.

The proof given in \cite{J1} for \eqref{J_sum_less_length} relied (quite heavily) on complex analysis.
In \cite{Ok} Okikiolu extended \eqref{J_sum_less_length} to $\Gamma \subset\R^d$ replacing complex analysis with  Euclidean geometry and some $\ell_2$ type computations.
The constants that hide behind the use of the symbol $\lesssim$ are exponential in $d$.  This
arises from the fact that a multi-scale dyadic grid is used, from some
accounting methods (which can in turn be related to the dyadic grid as 
well) and from ideas such as covering the unit sphere in $\R^d$ with balls of 
radius $\delta$.

Note that \eqref{J_sum_less_length} can be reformulated without defining $\D$ by 
\begin{gather}\label{J_int_less_length}
\int_0^\infty\int_\Gamma{\beta_\Gamma^2(\ball(x,At))\over\cH^1(\Gamma\cap \ball(x,t))}dx{dt}\lesssim
\cH^1(\Gamma) 
\end{gather}
where $A$ is a constant.  
See Lemma \ref{discretization} for some further details (but not all of them as we consider a different multiresolution in that lemma).

\subsection{New Results}\label{sect_new_results}
We prove a Hilbert space version of the above results.

Let a set $K\subset H$ be given.  For $n>n_0$ define $X^K_n \subset K$ to be a $2^{-n}$ net such that $X^K_n \subset X^K_{n+1}$.
Define a replacement for $\D$:
\begin{gather*}
\hat{\G}^K=\{Q=\ball(x,A2^{-n}):x\in  X^K_n, n \text{ an integer}, n\geq n_0\}.
\end{gather*}   
where $A>1$ is a constant and , $n_0$ is a (possibly negative) integer. 

We show (in Section \ref{sum_leq_thms})
\begin{theorem}\label{sum_beta_less_length}
\begin{gather*}
\suml_{Q\in\hat{\G}^K} \beta_\Gamma^2(Q)\diam(Q)
\lesssim \cH^1(\Gamma)
\end{gather*}
for any  connected set $\Gamma$ containing  $K$.
The constant behind the symbol $\lesssim$ depends only on the choice of $A$ (which can be given any value greater then $1$). In particular, the constant is independent of our choice of $\{X^K_n\}_{n\geq n_0}$ and the choice of $n_0$.
\end{theorem}
Equivalently (see Corollary \ref{cor_discretization}),
\begin{theorem}\label{int_beta_less_length}
\begin{gather*}
\int_0^\infty\int_\Gamma{\beta_\Gamma^2(\ball(x,At))\over\cH^1(\Gamma\cap \ball(x,t))}dx{dt}
\lesssim \cH^1(\Gamma)
\end{gather*}
for any  connected set $\Gamma$.
The constant behind the symbol $\lesssim$ depends only on the choice of $A$ (which can be given any value greater then $1$).
\end{theorem}

\begin{rem}
Our proof actually gives more information.  
Consider a set $E$ which is a  countable union of connected sets $\Gamma_i$
\begin{gather*}
E=\cup \Gamma_i.
\end{gather*}
Let $K\subset E$, and
\begin{gather*}
\G^{E,K} = \{Q \in \hat{\G^K} : \forall i, \quad \Gamma_i \cap (H\smallsetminus 4Q) \neq \emptyset\}.
\end{gather*}
Then inspection of the proof we give shows 
\begin{gather*}
\suml_{Q\in\G^{E,K}} \beta_E^2(Q)\diam(Q)
\lesssim \sum\cH^1(\Gamma_i)
\end{gather*} 
\end{rem}

\begin{rem}
For the statement and proof of theorem \ref{sum_beta_less_length} and for the previous remark we do not actually need the condition $X_{n+1}\supset X_n$ in the definition of $\hat{\G}^K$.  This condition is however used for the statement and proof of the following results.  
\end{rem}

In Section \ref{constructions}  we modify Jones' construction to give:
\begin{theorem}\label{construction_thm}
There is a constant $A_0$, such that for all $A>A_0$, 
for any set $K\subset H$ there exists
a  connected set $\Gamma_0 \supset K$ satisfying 
\begin{gather}\label{length_less_sum_beta}
\cH^1(\Gamma_0) \lesssim \diam(K) + \suml_{\hat{\G}^K} \beta_K^2(Q)\diam(Q).
\end{gather}
The constant behind the symbol $\lesssim$ depends only on the choice of $A$. In particular, the constant is independent of our choice of $\{X^K_n\}_{n\geq n_0}$.  We require $2^{-n_0}\geq \diam(K)$.
\end{theorem}
\begin{rem}
We would like to note recent independent work done by Immo Hahlomaa (see \cite{Ha}) containing  a generalization of Theorem \ref{construction_thm} to the setting of Metric Spaces! (Where one needs to use Menger curvature to define $\beta$.)  There is also work by Ferrari, Franchi and Pajot for a version of this theorem in certain geodesic spaces (such as the Heisenberg group) \cite{FFP}.
Analogs of Theorem \ref{sum_beta_less_length} and Theorem \ref{int_beta_less_length} can be obtained for Ahlfors-regular metric spaces \cite{RS-metric, Ha-2}. 
See the survey  \cite{my-TSP-survey} for some more details on the above results (without more than a hint of the proofs).
\end{rem}

As immediate corollaries (by combining the above theorems), we get the following results, which in $\R^2$ were the motivation for \cite{J1}.  
\begin{cor}
Let $\Gamma_0$ be as constructed in Theorem \ref{construction_thm}.
For $A>A_0$ and $n_0\leq -\log(\diam(K))$
\begin{gather*}
\cH^1(\Gamma_0) \lesssim \cH^1(\Gamma_{MST}).
\end{gather*}
\begin{gather*}
\diam(K)+ \suml_{Q\in\hat{\G}^K} \beta_K^2(Q)\diam(Q)
\sim \cH^1(\Gamma_{MST})
\end{gather*}
for any   set $K$, where $\Gamma_{MST}$ is a shortest connected set containing $K$.
\end{cor}
We show  the existence of $\Gamma_{MST}$ in Appendix \ref{appendix_MST}.
\begin{cor}\label{int_beta_sim_length}
For $A>A_0$
\begin{gather*}
\diam(\Gamma)+ \int_0^\infty\int_\Gamma{\beta_\Gamma^2(\ball(x,At))\over\cH^1(\Gamma\cap \ball(x,t))}dx{dt}
\sim \cH^1(\Gamma).
\end{gather*}
\begin{gather*}
\diam(\Gamma)+ \suml_{Q\in\hat{\G}^\Gamma} \beta_\Gamma^2(Q)\diam(Q)
\sim \cH^1(\Gamma)
\end{gather*}
for any  connected set $\Gamma$.
\end{cor}

This automatically gives that the relevant constants  in \cite{J1,Ok} need not be exponential in $d$, if $\D$ is replaced by $\hat{\G}$ (super-indexed correctly).  
The construction of $\Gamma_0$ remains the same as Jones', except for
one part (specifically, the case when $\beta > \epsilon$ becomes slightly more complicated).  
The proof of Theorem \ref{sum_beta_less_length} is done by a modification (as described in the following paragraphs) of Okikiolu's method which results in dimension free estimates.

One should note that in the case of one dimensional Uniformly Rectifiable sets 
(see \cite{DS} for definition of Uniformly Rectifiable) 
these results are obtained with much less difficulty by combining  \cite{Da,DS,J2,J1,Ok}.
The key idea is that using \cite{J2} one gets  that Ahlfors regular curves contain what is called `big pieces of chord-arc curves' (see \cite{DS} for a definition).
For chord-arc curves we have (using a modification of \cite{Ok})  desired estimates, which can be used with machinery from \cite{DS} to extend to Ahlfors regular curves. 
All this requires
an inspection of some proofs given in the above references, which results in the observation that, even though they are not stated as such, they are dimension independent for the relevant cases (or can be  made so with very minor modifications; for example \cite{Ok} can be made dimension independent in the case of chord-arc curves).  Inspecting the results in \cite{DS} was suggested to the author by Guy David.\\

\subsubsection*{Outline}
We prove  Theorems \ref{sum_beta_less_length} and \ref{int_beta_less_length} in Section \ref{sum_leq_thms}.   
We do this by considering the geometry of the set $\Gamma$ inside the different balls.
Let us give a vague intuitive description.\\

We call an arc $\tau$ (delimited by a given ball $Q$ and contained in $\Gamma\cap Q$) an `almost flat arc' if  $\beta(\tau)$ is small in comparison with $\beta(Q)$ (i.e. $\tau$ is close to a straight line segment). 
For a given $Q$ the collection of these arcs is called $S_Q$.  We also designate an arc going through the center of $Q$ (existing by the definition of $\hat{\G}$) by $\gamma_Q$.

\begin{figure}[h]
\begin{center}
$\begin{array}{ccc}
\scalebox{0.1}{\includegraphics*[0in,0in][8in,8in]{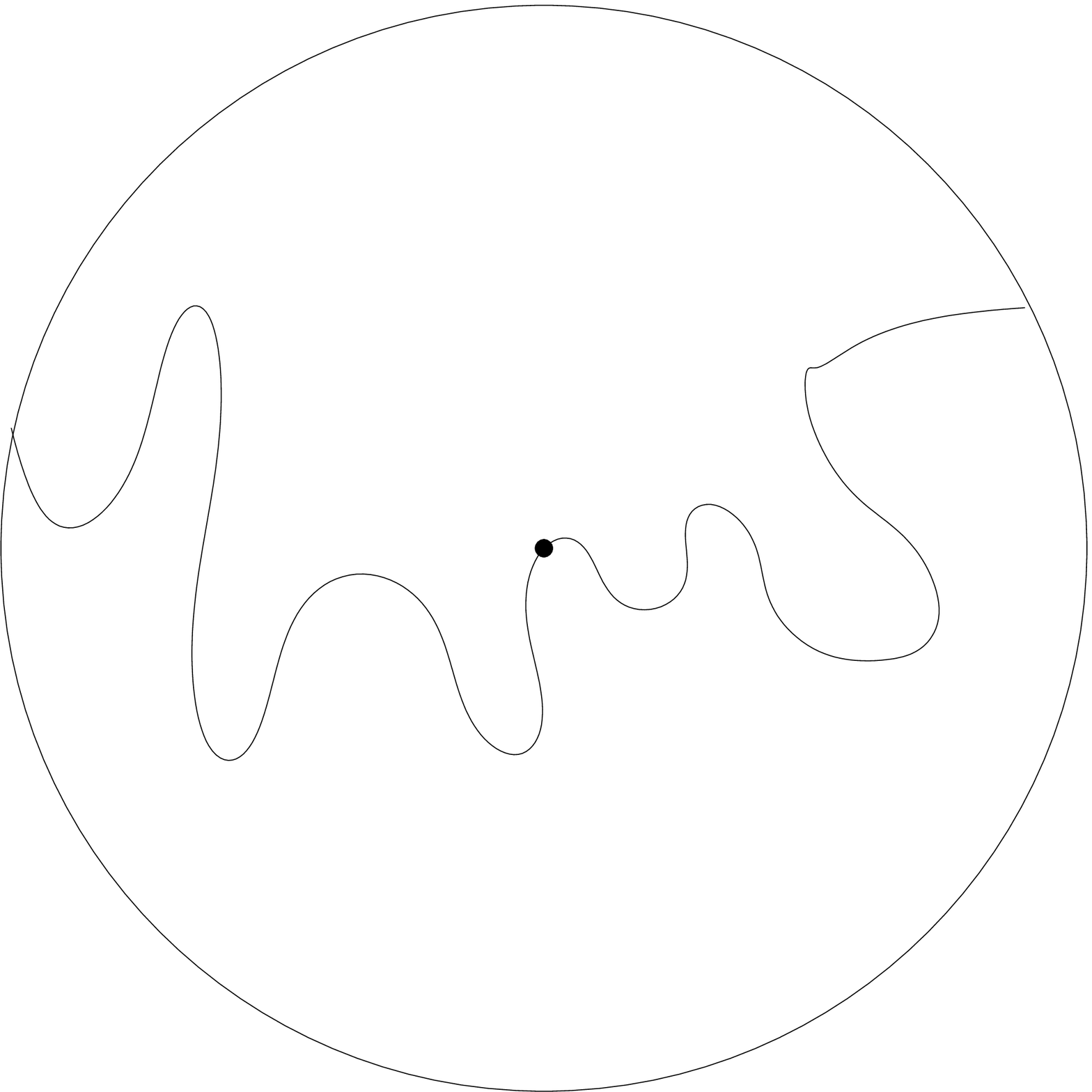}}&
\scalebox{0.1}{\includegraphics*[0in,0in][8in,8in]{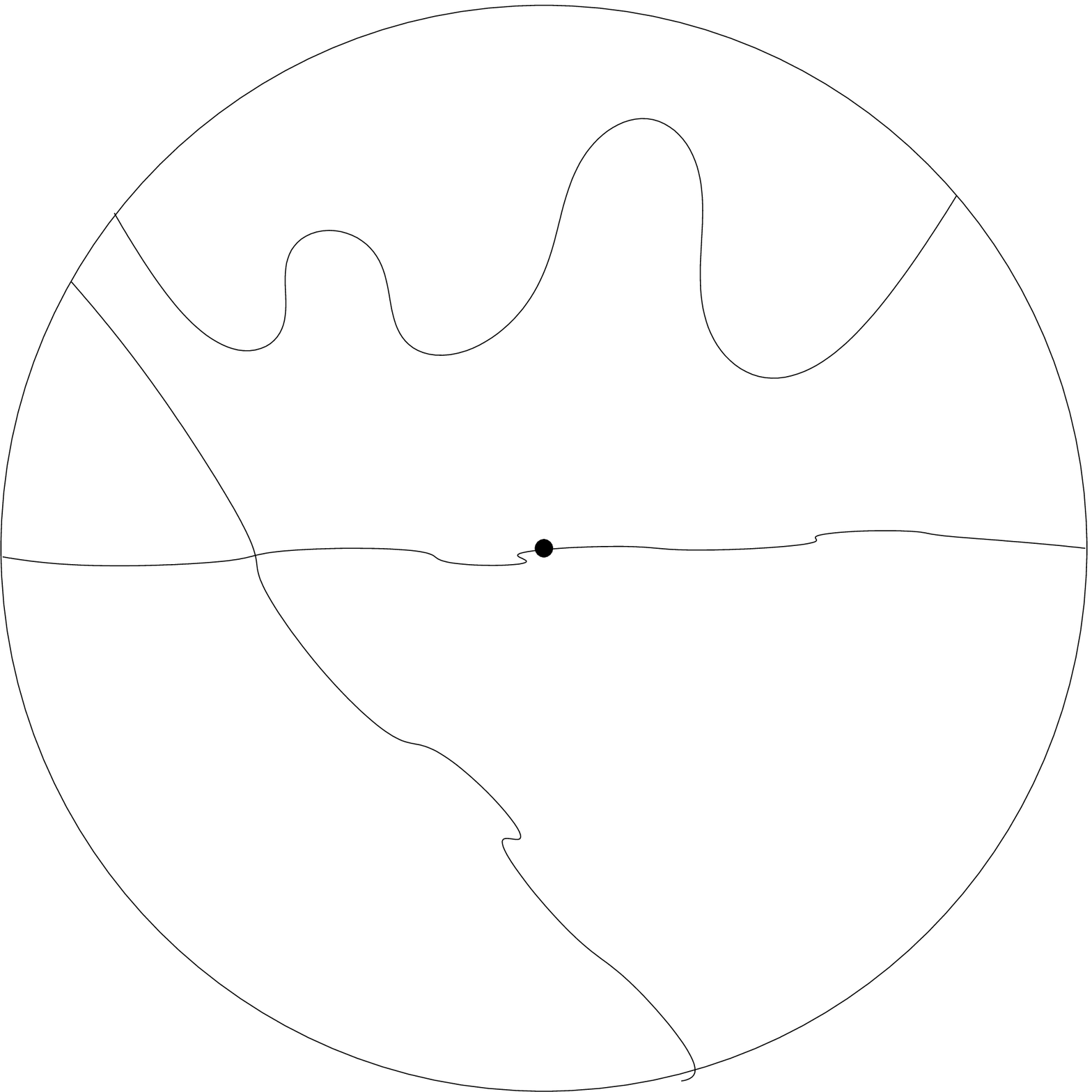}}&
\scalebox{0.1}{\includegraphics*[0in,0in][8in,8in]{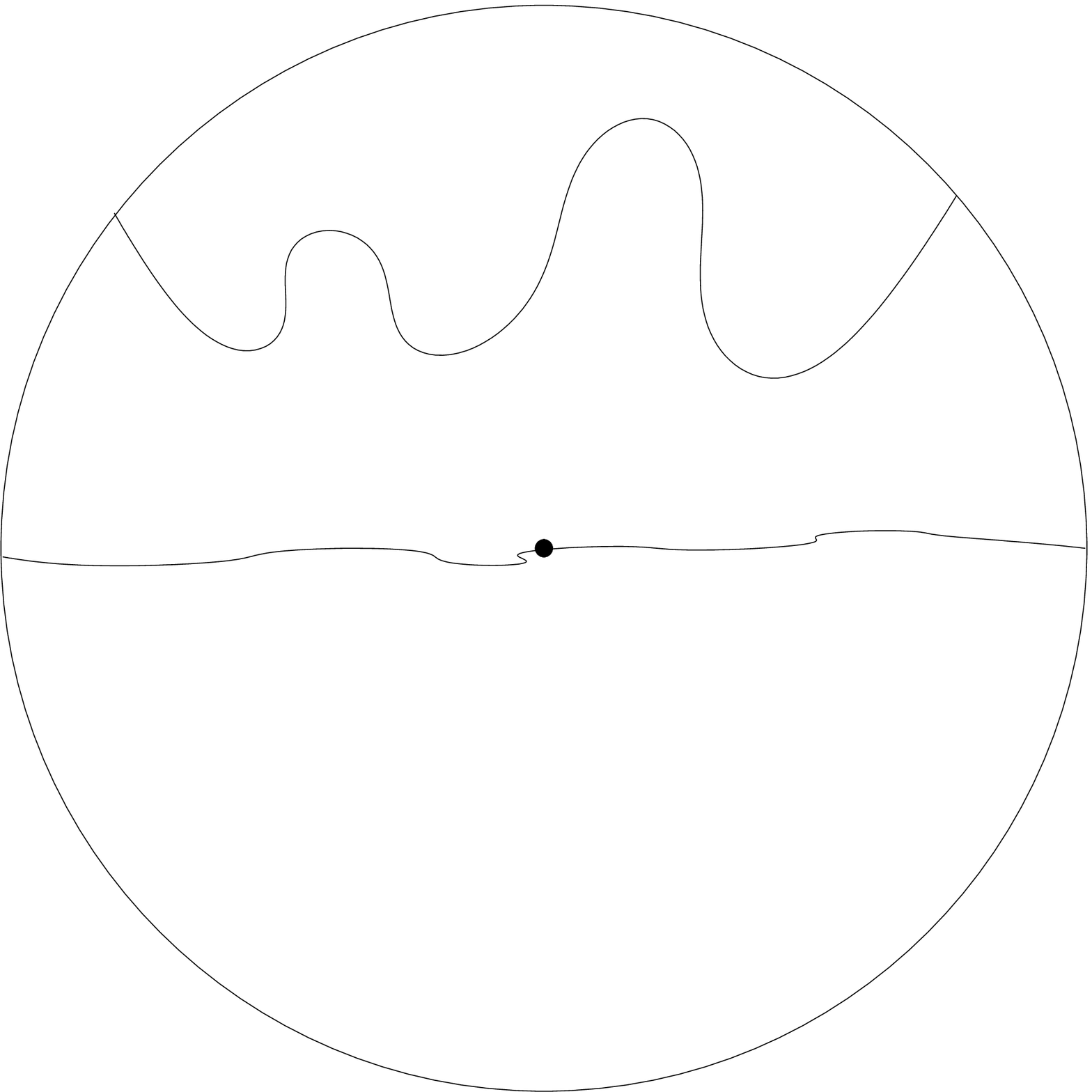}}\\
(a)&(b)&(c)
\end{array}$
\end{center}
\caption{Three examples of  balls}
\label{mini_A_examples}
\end{figure}

In subsection \ref{curvy_arcs} we discuss balls $Q$ for which either $\gamma_Q$ is not an `almost flat arc'  (`non-flat arc' ) or $\beta(Q)$ is not controlled by $\beta_{S_Q}(Q)$ ($=\beta$ restricted to  the `almost flat arcs'). 
The latter  balls also contain an arc $\tau$ which is `very non-flat'  in the sense that it is `non-flat' enough so that $\beta(\tau)$  controls $\beta(Q)$  
(see Figure \ref{mini_A_examples} (a) and (c)). 
We make use of this by employing ideas of Okikiolu, as well as ideas similar to ones of  G. David and M. Christ (see e.g. \cite{Ch}, and \cite{Da} page 93 for a simple version).  Note however that one main difference with the case of Christ and David is the lack of  homogeneity assumption on $\Gamma$! 
This makes Okikiolu's ideas harder to use. Subsection \ref{curvy_arcs} corresponds to the first half of \cite{Ok}. 

In subsection \ref{almost_straight_lines} we discuss the rest of the balls.  These balls $Q$ satisfy  $\beta(Q)  \lesssim\beta_{S_Q}(Q)$ (see Figure \ref{mini_A_examples} (b)).
The key idea here is to use the curve itself as the notebook for the bookkeeping (as in \cite{Ok,Le,BJ}).
This is done explicitly in this subsection (whereas in this  preceding subsections this idea is used implicitly). 
This subsection corresponds to  the second part of \cite{Ok}, where Okikiolu allots segments whose length controls $\beta(Q)\diam(Q)$.  
We substitute allotting segments by allotting  densities.  
We do so by constructing for each ball $Q \in \hat{\G}$, a weight $w_Q$ supported in $Q$, satisfying
$\int_Q w_Qd\length \geq \beta(Q)\diam(Q)$ and  $\suml_Q w_Q(x) \leq C$ for almost every $x\in \Gamma$.
An important point is that the construction of each $w_Q$ is done in a multiscale fashion (as a martingale), which allows  the assurance of the above properties.  This assurance is  not so straight forward and most of subsection \ref{almost_straight_lines} is devoted to it.
This gives us that for these balls $\sum\beta_\Gamma(Q)\diam(Q) \lesssim \cH^1 (\Gamma)$ (note that we do not square the terms!).  \\

We prove  Theorem \ref{construction_thm} in Section \ref{constructions}.  We do so by a construction which is a modification of the {\it farthest insertion} algorithm.  This is not far from what is written in \cite{J1}.     

\subsubsection*{Constants and a computational note}
We will fix certain constants in the following proofs.  Some of them will depend on each other.  In particular, in subsection \ref{farthest_local_hilbert} we will obtain a value $A_0$ so that we will require $A>A_0$. $A_0=200$ will suffice. Then we will use a (any) choice of $A$ and derive from it a choice for the constant that we name $\epsilon_2$, which first appears in subsection \ref{preliminaries}.  In subsection \ref{almost_straight_lines} we have a constant $C$ whose choice depends on $A$.  The constant $J$ appears several times when we wish to skip (Jump over) scales.  Only in subsection  \ref{almost_straight_lines} is it required that it depend on $A$.  

All other named constants are independent of $A$.  

The constant $n_0$ introduced in this section is not used in any other constant! It is used solely as a starting point for an inductive argument. 
 
By following the proofs one gets that  
the dependence of the constant hiding behind the symbol $\sim$ in Corollary \ref{int_beta_sim_length}  on $A$ is 
$A^{9\over 2}\log{A}$.  We have made no effort to get this dependence to be as minimal as possible while proving the theorems.  The reason this is of interest is that one may find this useful when trying to use these theorems in an $\R^d$ numerical setting, and hence go back to considering dyadic cubes (which renders the choice of  $A\sim\sqrt{d}$ reasonable).
%
%
%
%
%
%
%
%
%
%
%
%
%
%
%
%
%
%
%
%
%
%
%
%
%
%
%
%

\section{Acknowledgements}
The work presented here comes from the content of the author's PhD dissertation, which was done under the direction of Peter  Jones at Yale University.
The author would like to thank Peter Jones, for the many hours of discussion.  The author would also like to thank Guy David and Christopher Bishop for their many useful comments on the PhD dissertation.
Finally, the author would like to thank the referee for many kind and helpful comments.\\
\section{Proof of Theorems \ref{sum_beta_less_length} and \ref{int_beta_less_length}}\label{sum_leq_thms}
%
%
%
\subsection{Preliminaries, Notation and Definitions}\label{preliminaries}
\begin{rem}
We would like to show
\begin{gather}\label{09022005}
\sum\limits_{Q \in \hat{\G}}\beta(Q)^2\diam(Q) \leq C \cH^1 (\Gamma).
\end{gather}
We have
\begin{gather*}
\lim\limits_{n \to \infty}\sum\limits_{Q \in \hat{\G} \atop center(Q) \in X_n}\beta_{X_n}(Q)^2\diam(Q)=
\sum\limits_{Q \in \hat{\G}}\beta(Q)^2\diam(Q).
\end{gather*}
If we choose $d=\sharp(X_n)$ then we can project the problem to $\R^d$, and so if we prove \eqref{09022005}
for $\Gamma \subset \R^d$ with $C$ independent of $d$ we are done.
Hence one should note that we may just as well in subsections \ref{preliminaries} - \ref{almost_straight_lines}  assume that we are working in $\R^d$, not in $H$.  This will however, be of no consequence to us, as our proof works in $H$.
\end{rem}

Assume $K \subset\Gamma\subset H$ as in the statement of the theorems. 
$A$ is fixed to a constant larger or equal to the constant $A_0$ that we get from Chapter \ref{constructions}.  We will omit the superscripts/subscripts $K,\Gamma$ whenever possible to simplify notation.\\

We start with a discretization lemma.
\begin{lemma}\label{discretization}
Assume that $K\subset \Gamma$ and $\cH^1(B(x,t)\cap {\Gamma}) >0, \forall x\in {\Gamma}, t >0$.
Let $8A'\leq A \leq {1\over 8}A''$.
Then
\begin{gather*}
(i) \int_0^\infty\int_{\Gamma}{\beta_\Gamma^2(A'\ball(x,t))\over\cH^1({\Gamma}\cap \ball(x,t))}dx{dt} \lesssim
\sum\limits_{Q \in \hat{\G^\Gamma}}\beta_\Gamma(Q)^2\diam(Q)\\
(ii) \int_0^\infty\int_{\Gamma}{\beta_\Gamma^2(A''\ball(x,t))\over\cH^1({\Gamma}\cap \ball(x,t))}dx{dt} \gtrsim
\sum\limits_{Q \in \hat{\G^K}}\beta_\Gamma(Q)^2\diam(Q).
\end{gather*}
\end{lemma}
\begin{proof}
Notice that for a ball $B$ we have $\beta(B)\leq C\beta(CB),\forall C\geq1$.  
\begin{eqnarray*}
\int_0^\infty\int_{\Gamma}{\beta^2(A'\ball(x,t))\over\cH^1({\Gamma}\cap \ball(x,t))}dx{dt} &\lesssim&
\int_0^\infty\int_{\Gamma}{\beta^2(A'\ball(x,2t))\over\cH^1({\Gamma}\cap \ball(x,2t))}dx{dt}\\
 &\lesssim&
\suml_{n\in\Z}2^{-n}\int_{\Gamma}{\beta^2(4A'\ball(x,2^{-n}))\over\cH^1({\Gamma}\cap 2\ball(x,2^{-n}))}dx\\
&\lesssim&
\suml_{n\in\Z}2^{-n} \suml_{x\in X_n^\Gamma} \beta^2(8A'\ball(x,2^{-n}))
	{\cH^1({\Gamma} \cap \ball(x,2^{-n}))\over\cH^1({\Gamma} \cap \ball(x,2^{-n}))}\\
&\lesssim&
\suml_{Q \in \hat{\G^\Gamma}}\beta(Q)^2\diam(Q)
\end{eqnarray*}
as long as $A\geq 8A'$. The change of variable $t\to\half t$ was used for the first inequality.   
Conversely,
\begin{eqnarray*}
\suml_{Q \in \hat{\G^K}}\beta(Q)^2\diam(Q)
&=&
\suml_{n\in\Z}2A2^{-n} \suml_{x\in X_n^K} \beta^2(\ball(x,A2^{-n}))\\
&\lesssim&
\suml_{n\in\Z}2^{-n} \int_{x\in {\Gamma}} \beta^2(2\ball(x,A2^{-n}))
	{1\over\cH^1({\Gamma} \cap \half \ball(x,2^{-n}))}dx\\
&\lesssim&
\int_0^\infty \int_{x\in {\Gamma}} \beta^2(4\ball(x,At))
	{1\over\cH^1({\Gamma} \cap \half \ball(x,t))}dxdt\\
&\lesssim&
\int_0^\infty \int_{x\in {\Gamma}} {\beta^2(8A\ball(x,t))
	\over\cH^1({\Gamma} \cap \ball(x,t))}dxdt
\end{eqnarray*}
so take $A''\geq8A$. The change of variable $t\to2t$ was used for the last inequality.
\end{proof}
This immediately gives us
\begin{cor}\label{cor_discretization}
Theorem \ref{int_beta_less_length} and Theorem \ref{sum_beta_less_length} are equivalent.
\end{cor}

For the remainder of this section we will  concern ourselves only with Theorem \ref{sum_beta_less_length}.
We state some point-set topology lemmas. For completeness  we give proofs for these lemmas in the appendix.
\begin{lemma}\label{closure-length}
Assume $\Gamma$ is connected.  Then $\cH^1(\Gamma)=\cH^1(\Gamma^{closure})$.
\end{lemma}
\begin{lemma}\label{finite_length_then_cpt}
Assume $\Gamma\subset H$ is a closed connected set with $\cH^1(\Gamma)< \infty$.
Then $\Gamma$ is compact.
\end{lemma}
\begin{lemma}\label{cpt}
Let $C_1,C_2 >0$ be given.
Given a compact connected set $\Gamma \subset H$ the set
$E:=\{x\in H:x=tx_1 + (1-t)x_2, x_i \in \Gamma, -C_1\leq t \leq C_2\}$
is compact.
\end{lemma}
\begin{lemma}\label{parametrization}
Let $\Gamma \subset H$ be a compact connected  set of finite length.  
Then we have a Lipschitz function $\gamma:[0,1] \to H \st Image (\gamma)=\Gamma$ 
and $\norm{\gamma}_{Lip} \leq 32\cH^1 (\Gamma)$
\end{lemma}
\begin{cor}\label{cor-parametrization}
Let $\Gamma \subset H$ be a compact connected  set of finite length.  
Then we have a Lipschitz function $\gamma:\T \to H \st Image (\gamma)=\Gamma$ 
and $\norm{\gamma}_{Lip} \leq 32\cH^1 (\Gamma)$
\end{cor}

\noindent
Proofs for the above lemmas can be found in the appendix.
\medskip

Since Theorem \ref{sum_beta_less_length}  is  trivially satisfied for $\Gamma$ satisfying $\cH^1(\Gamma)=\infty$ we will assume $\cH^1(\Gamma)< \infty$.  
We may replace $\Gamma$ by its closure without loss of generality for the purpose of proving this theorem.  
This will not affect  the Jones-$\beta$ numbers, the connectedness of $\Gamma$, or the  length of $\Gamma$.  
Hence we will assume  $\Gamma$ is compact from now on.
By re-scaling we may also assume $\diam(\Gamma)\leq1$ and $n_0=0$.

\smallskip

Using Corollary \ref{cor-parametrization}, we fix a parameterization $\gamma$ for $\Gamma$, $\gamma:\mathbb{T} \longrightarrow H$,
such that we have $\length(\gamma) \leq 32\cH^1(\Gamma)$.  {\bf From here on we also use $\length(\cdot)$ as the push-forward by $\gamma$ of the arc-length measure}.

We will show (for a given $A$) that Theorem \ref{sum_beta_less_length} holds, i.e. 
\begin{gather*}
\sum\limits_{Q \in \hat{\G}}\beta(Q)^2\diam(Q) \leq C \cH^1 (\Gamma).
\end{gather*}
We recall that  we have (after re-scaling)
\begin{gather}\label{g-def}
\hat{\G} = \{Q =\ball(x,A2^{-n}),x \in X_n; n \geq 0\}.
\end{gather}
We define
\begin{gather*}
\G = \{Q \in \hat{\G} :\Gamma \cap (H\smallsetminus 4Q) \neq \emptyset\}.
\end{gather*}
($\G$ is the collection of balls $Q$ that are small enough so that $\Gamma$ must exit $4Q$.)

Consider $\hat{\G} \smallsetminus \G$.
\begin{lemma}
$\sum\limits_{Q \in \hat{\G} \smallsetminus \G}\beta(Q)^2\diam(Q) \leq C \cH^1 (\Gamma)$.
\end{lemma}
\begin{proof}
Set $L=\length(\gamma)$ 
(the arc-length of a parameterization $\gamma$ of $\Gamma$ assured by Corollary 
\ref{cor-parametrization})
and $D=\diam(\Gamma)$.  
We have at most $\frac{8\cdot2AL}{D}$ balls of diameter $\eighth D$ in $\hat{\G} \smallsetminus \G$,
as the centers of these balls are at least $\frac{D}{8\cdot2A}$ apart and are along a curve of length $L$.
Similarly, we have at most:\\
\begin{eqnarray*}
\frac{4\cdot2AL}{D}&\mbox { balls of diameter }\quarter D &\mbox { in }\hat{\G} \smallsetminus \G.
\mbox { They are of }\beta \leq 1.\\
\frac{2\cdot2AL}{D} &\mbox { balls of diameter }\half D &\mbox {  in }\hat{\G} \smallsetminus \G.
\mbox { They are of } \beta \leq 1.\\
\frac{2AL}{D} &\mbox { balls of diameter } D &\mbox { in }\hat{\G} \smallsetminus \G.
\mbox { They are of }\beta \leq 1.\\
\frac{2AL}{2D}&\mbox { balls of diameter } 2D &\mbox {  in }\hat{\G} \smallsetminus \G.
\mbox { They are of }\beta \leq \half.\\
\frac{2AL}{4D} &\mbox { balls of diameter } 4D &\mbox { in }\hat{\G} \smallsetminus \G.
\mbox { They are of }\beta \leq \quarter.\\
\vdots\\
\frac{2AL}{AL} &\mbox { balls of diameter } AL &\mbox {  in } \hat{\G} \smallsetminus \G.
\mbox { They are of }\beta \leq \frac{D}{AL}\\
1 &\mbox { ball of diameter } 2AL &\mbox {  in } \hat{\G} \smallsetminus \G.
\mbox { It is  of }\beta \leq \frac{D}{2AL}\\
1 &\mbox { ball of diameter } 4AL &\mbox {  in } \hat{\G} \smallsetminus \G.
\mbox { It is  of }\beta \leq \frac{D}{4AL}\\
\vdots
\end{eqnarray*}
(Since the centers must be  along $\gamma$, which is of length $L$). 
Hence 
\begin{gather*}
\suml_{Q\in \hat{\G} \smallsetminus \G} \beta(Q)^2\diam(Q) \leq 
\suml_{n=-3}^{\log(\frac{AL}{D})} AL 2^{-2n}  + \sum_{\log(\frac{AL}{D})}^\infty AL 2^{-n}
\leq 32AL + D \leq (32A+1)L.
\end{gather*}
\end{proof}

We need some more notation:
\begin{eqnarray}\label{lambda-definition}
\Lambda(Q)&:=&
	\{\tau=\gamma|_{[a,b]}: [a,b]\subset \T; [a,b] 
		\text { a connected component of }\gamma^{-1}(\Gamma \cap Q)\}.
\end{eqnarray}
We will freely use $\tau\in \Lambda(Q)$ as both  a parameterization of an arc (given by restriction of $\gamma$), and its image.
We will denote by $\diam(\tau)$ the diameter of the image of $\tau$.

We define for an arc $\tau:[a,b]\to H$
\begin{eqnarray}\label{bt-definition}
\bt(\tau):=\supl_{t \in [a,b]} \frac{\dist(\tau(t),[\tau(a),\tau(b)])}{\diam(\tau)},
\end{eqnarray}
where $[x,y]$ is the straight line segment connecting $x$ and $y$.
(This is how we define the Jones $\beta$ number of an arc).

\begin{figure}[p]
\begin{center}
$\begin{array}{c}
\scalebox{0.35}{\includegraphics*[0in,0in][8in,8in]{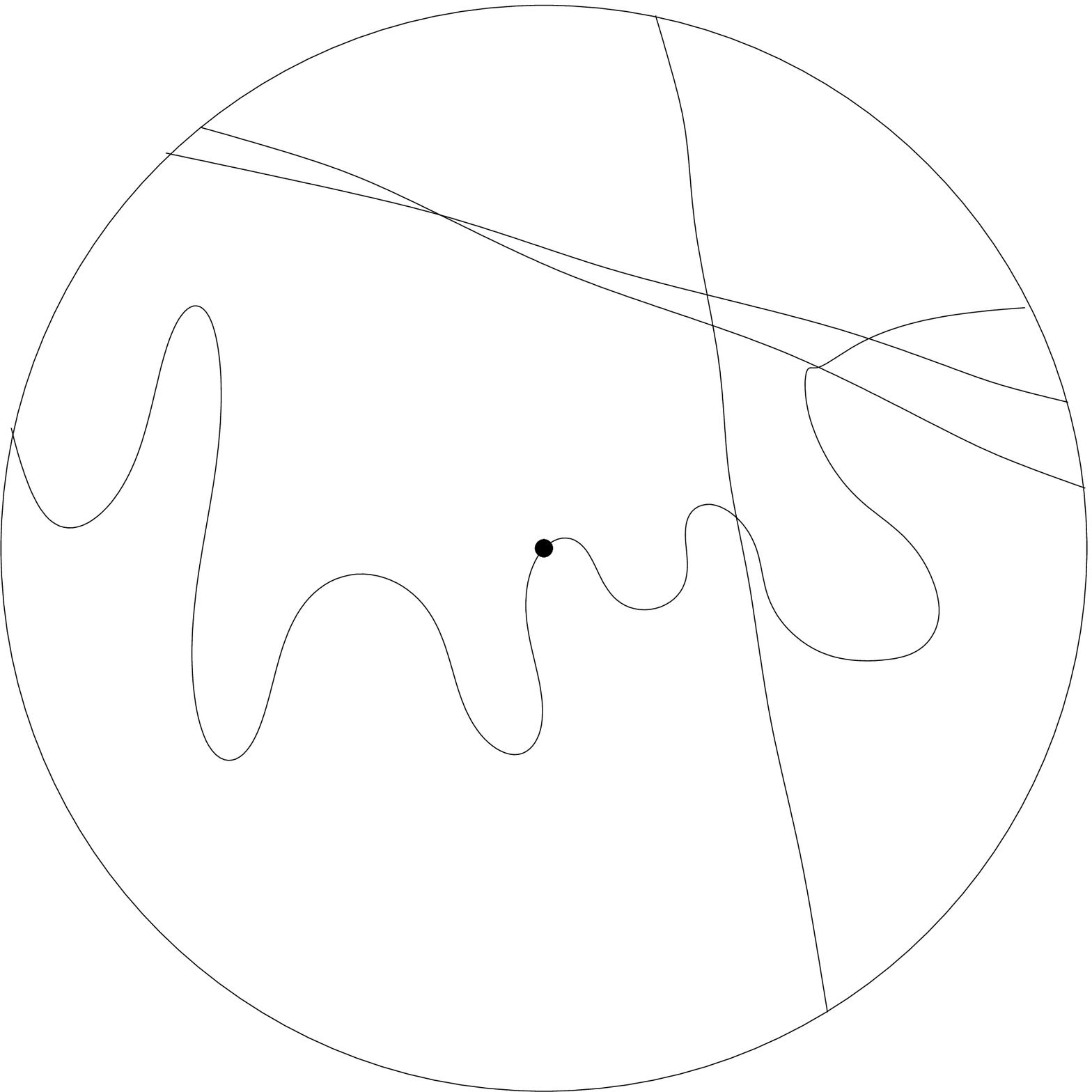}}\\
\scalebox{0.35}{\includegraphics*[0in,0in][8in,8in]{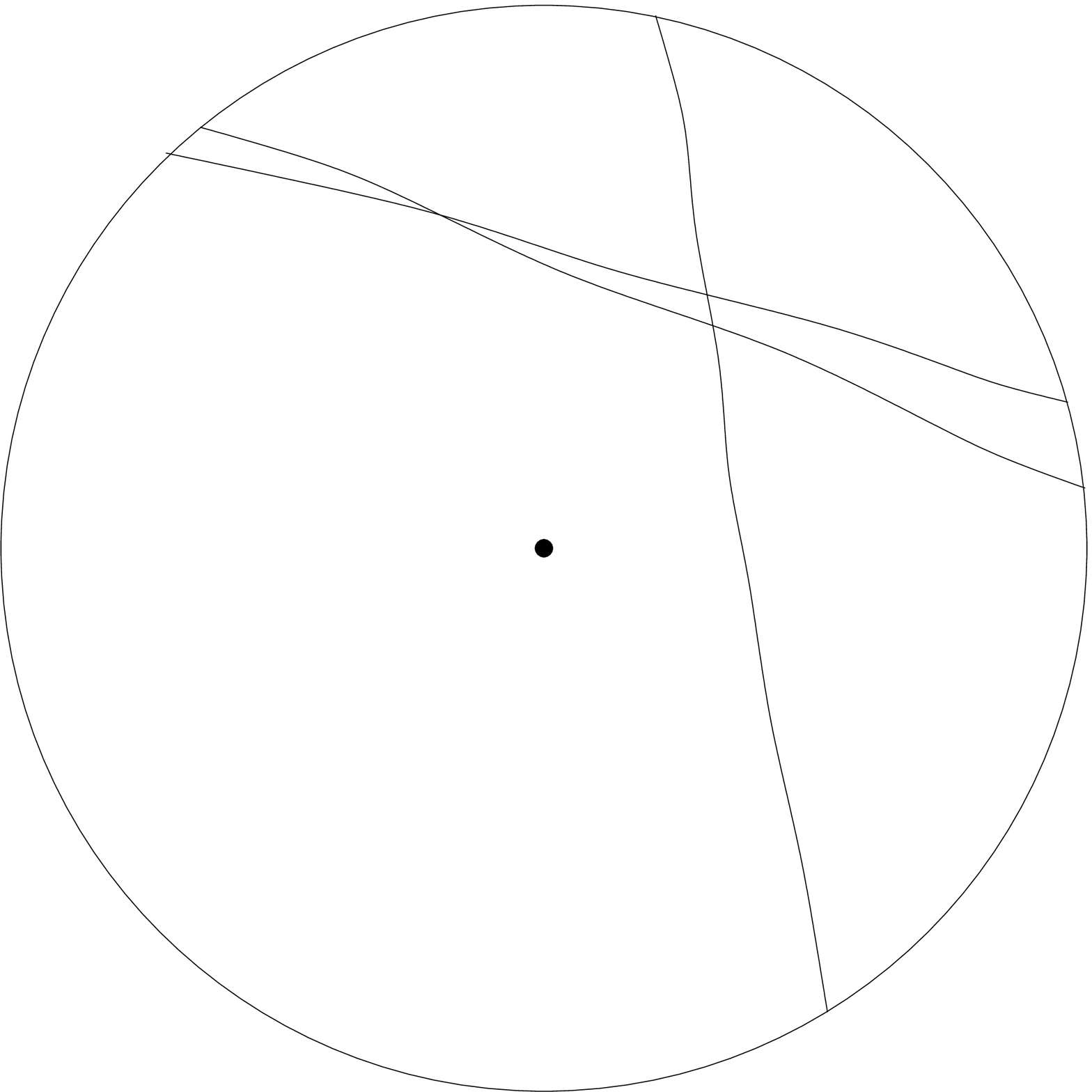}}
\end{array}$
\end{center}
\caption{Example of $Q$ (top) and $S_Q$ (bottom)}
\label{S_example}
\end{figure}

Consider $\tau\in \Lambda(Q)$.  
We call $\tau$ \textit{almost flat} iff
\begin{eqnarray*}
\bt(\tau)\leq \epsilon_2\beta(Q)\,,
\end{eqnarray*}
where $\epsilon_2$  is a constant which will be fixed  in subsection \ref{almost_straight_lines}, (with $A\epsilon_2$ independent of $A$ and sufficiently small).
Set:
\begin{eqnarray}\label{S_Q-definition}
S_Q:=\{\tau \in \Lambda(Q): \bt(\tau) \leq \epsilon_2 \beta(Q)\}.
\end{eqnarray}
This is the collection of \textit{almost flat} arcs  in $\Lambda(Q)$.  See Figure \ref{S_example}.

Consider all balls $Q \in \G \cup 2\G \cup 4\G$.  For each of them, we fix $\gamma_Q \in \Lambda_Q$ an arc containing the center of $Q$. 
If there is more then one option, choose so that if $Q_1=\half Q_2$ then 
$\gamma_{Q_1} \subset \gamma_{Q_2}$. This can be easily done by 
working top-down (as opposed to $\hat{\G}$, $\G$ has a coarsest scale). 

%
\begin{figure}[p]
\begin{center}
$\begin{array}{c}
\scalebox{0.25}{\includegraphics*[0in,0in][8in,8in]{pics/A1.eps}}
\put(-120,70){$\gamma_Q$}\\
\mbox{{An example of a $\G_1$ ball}}\\
\\\\
\scalebox{0.25}{\includegraphics*[0in,0in][8in,8in]{pics/A2.eps}}
\put(-70,70){$\gamma_Q$}\\
\mbox{{An example of a $\G_2$ ball}}\\
\\\\
\scalebox{0.25}{\includegraphics*[0in,0in][8in,8in]{pics/A3.eps}}
\put(-120,70){$\gamma_Q$}\\
\mbox{{An example of a $\G_3$ ball}}\\
\end{array}$
\end{center}
\caption{Examples of the three different types of balls}
\label{A_examples}
\end{figure}

\noindent
Set for $j\in\{0,1,2\}$ and  $\epsilon_1$, $\epsilon_2$ which will be fixed  in subsection \ref{almost_straight_lines} (with $\epsilon_1$  independent of $A$ and sufficiently small)
\begin{eqnarray*}
\G_1^j&=&\{Q\in \G: 
                                \bt(\gamma_{2^jQ})> \epsilon_2\beta(2^jQ)\}\\
\G_2^j&=&\{Q\in \G:
                                \bt(\gamma_{2^jQ})\leq \epsilon_2\beta(2^jQ); 
                                \beta_{S_{2^jQ}}(2^jQ) > \epsilon_1\beta(Q)\}\\
\G_3^j&=&\{Q\in \G:
                                \bt(\gamma_{2^jQ})\leq \epsilon_2\beta(2^jQ); 
                                \beta_{S_{2^jQ}}(2^jQ)  \leq \epsilon_1\beta(Q)\}
\end{eqnarray*}
where by $\beta_{S_{2^jQ}}(\cdot)$ we mean $\beta_{\cup \{\tau:\tau\in S_{2^jQ}\}}(\cdot)$. This is an abuse of notation we will keep on using throughout this essay.
Clearly for every $j$, $\G=\G_1^j \cup \G_2^j \cup \G_3^j$ 
and $\gamma_{2^jQ} \in S_{2^jQ}$, for all $Q\in \G_2^j \cup \G_3^j$.
See Figure \ref{A_examples}  for examples of the above sets. 
For most of this essay  the reader can be content with simply thinking of $j=0$ (the only exception will be $\G_3^j$ ).\\

We will show 
\begin{gather}\label{beta_over_A_1}
\sum\limits_{Q \in \G_1^j}\beta(Q)^2\diam(Q) \leq C \cH^1 (\Gamma) 
\end{gather} 
for any $j\in\{0,1,2\}$ in subsection \ref{curvy_arcs}.
We will also show 
\begin{gather}\label{beta_over_A_3}
\sum\limits_{Q \in \G_3^1\cap \G_3^2\cap\G_3^3}\beta(Q)^2\diam(Q) \leq C \cH^1 (\Gamma) 
\end{gather} 
in subsection \ref{curvy_arcs}.
We will show 
\begin{gather}\label{beta_over_A_2}
\sum\limits_{Q \in \G_2^j}\beta(Q)^2\diam(Q) \leq C \cH^1 (\Gamma) 
\end{gather} 
for any $j\in\{0,1,2\}$ in subsection \ref{almost_straight_lines}.

This will give us.
\begin{gather*}
\sum\limits_{Q \in \G}\beta(Q)^2\diam(Q) \leq C \cH^1 (\Gamma) 
\end{gather*}

%
%
\subsection{Non-Flat Arcs}\label{curvy_arcs}
In this subsection we prove \eqref{beta_over_A_1} and \eqref{beta_over_A_3}.
The tools developed in 
this subsection will also be used in subsection \ref{almost_straight_lines}.\\

\begin{rem}  
In the following we will be discussing various 
sub-arcs. They are parameterized by the global parameterization of $\Gamma$.
It is important that when we want to say something about the intersection 
or union of two sub-arcs, that we talk about their \textbf{domain} and
not their \textbf{range}!
In contrast, when we discuss the diameter of an arc, we will be discussing the diameter of its \textbf{image}.
\end{rem} 
   
Even though the setup is slightly different, the proof for the following lemma is copied almost word for word from \cite{Ok}.

\begin{lemma}\label{L2_lemma}
Suppose we are given a family of sub-arcs $\F=\bigcup\limits_{i=0}^\infty\F_i$  with the following properties:
\begin{verse}
(1) $\tau' \in \F_{n+1} \then \exists! \tau \in \F_{n} \st \tau' \subset \tau$\\
(2) $\tau \in \F_n \then 2^{-nJ} \leq \diam(\tau) \leq A2^{-nJ+2}$\\
(3) $\tau,\tau' \in \F_n \then \sharp(\tau \cap \tau') \in \{0,1,2\}$ 
(the intersection is an empty set, a single point, or two points)\\
(4) $\bigcup\limits_{\F_0} \tau = \bigcup\limits_{\F_n} \tau,$ $\forall n$
\end{verse}
(we will call such a family a {\bf filtration}).
Then we have:
\begin{gather*}
\sum\limits_{\tau \in \F} \bt(\tau)^2\diam(\tau) \lesssim \length(\bigcup\limits_{\F_0} \tau).
\end{gather*}
\end{lemma}
\begin{proof}
Set for $\tau \in \F_n$, 
\begin{gather*}
\F_{\tau,k} = \{\tau' \subset \tau: \tau' \in \F_{n+k}\}.
\end{gather*}
Set for $\tau \in \F_n \st \tau:[initial,final] \to H$, 
\begin{gather*}
I_{\tau} =[\tau(initial),\tau(final)]
\end{gather*} 
and  
\begin{gather*}
d_{\tau} = 
\supl_{\tau' \in \F_{\tau,1} \atop x \in I_{\tau'}} \dist(x,I_{\tau}).
\end{gather*} 
This is in fact a maximum and not a supremum:  $\bigcup\limits_{\tau' \in \F_{\tau,1}}I_{\tau'}$ is compact by the compactness of 
$\tau$.  
For each $\tau$ and $k$, let $\tau_k\in \F_{\tau,k}$ be chosen such that $d_{\tau_k}$ is maximal among all arcs in $\F_{\tau,k}$.  Again, this is a maximum and not a supremum.
%
%
\begin{figure}[p]
\begin{center}
\scalebox{.9}{\includegraphics*[0in,0in][5in,8in]{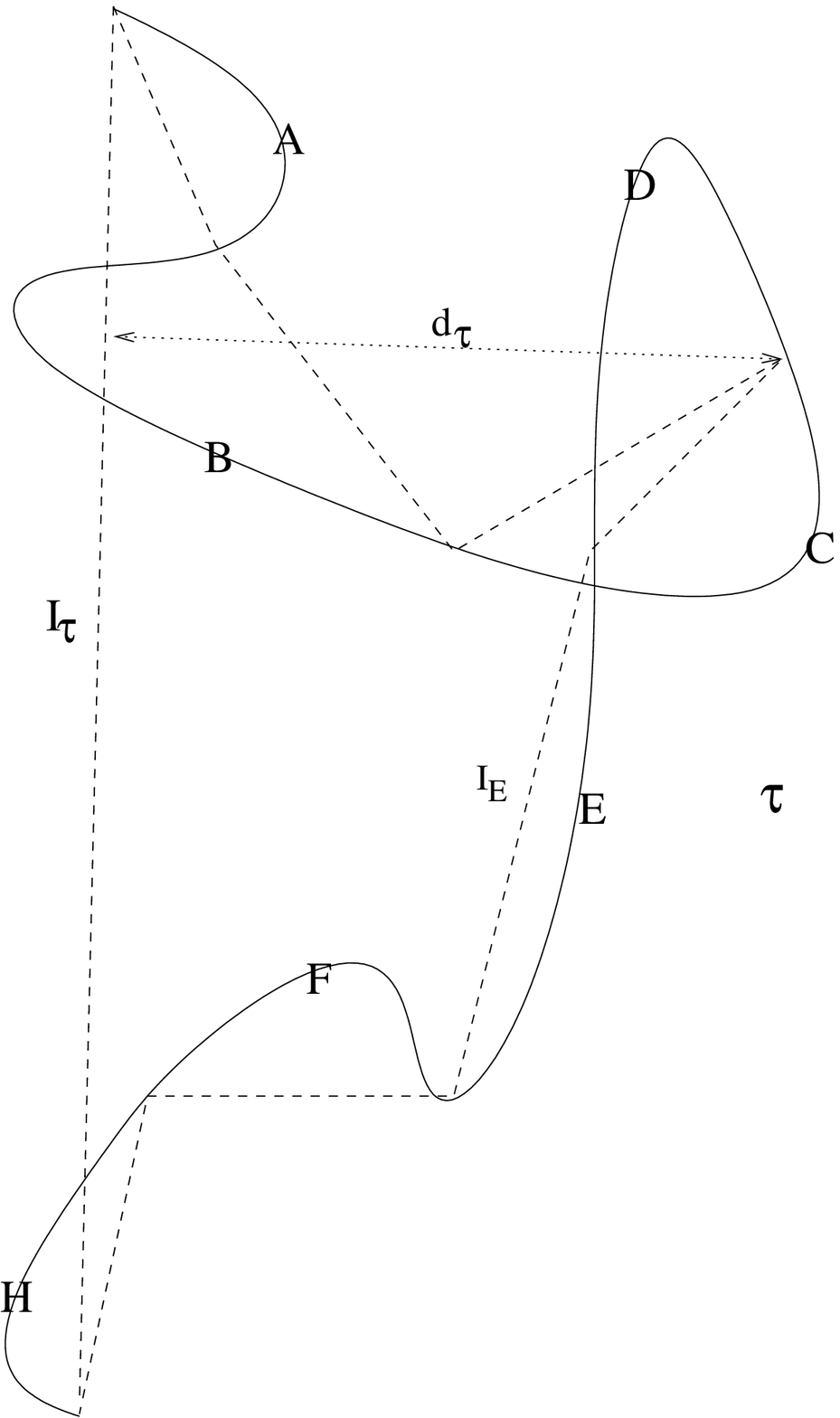}}
\end{center}
\caption[Picture for lemma \ref{L2_lemma}]{(Here $A,B,C,D,E,F,G,H \in F_{\tau,1}$ and 
$ A\cup B \cup C\cup D\cup F \cup G \cup H = \tau$)}
\label{ok_fig}
\end{figure}

We have  for $\tau \in \F_n$: 
\begin{gather}\label{temp}
\bt(\tau) \diam (\tau) \leq\sum\limits_{k=0}^\infty d_{\tau_k}
\end{gather} 
by the following.
Consider a sequence: $\tau^0=\tau$, $\tau^{k+1} \in \F_{\tau^k,1}$.   
We have 
\begin{gather*}
\bt(\tau)\diam(\tau)= \sum\limits_0^\infty (\bt(\tau^k)\diam(\tau^k)-\bt(\tau^{k+1})\diam(\tau^{k+1}))
\end{gather*}
by the fact that  we have a telescoping series with the summand going to 0 by (2).
If we choose $\tau^{k+1}\in \F_{\tau^k,1}$ such that $\bt(\tau^{k+1})\diam(\tau^{k+1})$ is maximal among all arcs in $\F_{\tau^k,1}$, then
\begin{gather*}
d_{\tau_k} \geq (\bt(\tau^k)\diam(\tau^k)-\bt(\tau^{k+1})\diam(\tau^{k+1}))
\end{gather*} 
by the triangle inequality and our choices of $\tau_k$ and $\tau^{k+1}$.
Thus we have \eqref{temp} 

We also have 
\begin{gather}\label{temp2}
\frac{d_\tau^2}{\diam(\tau)} 
\lesssim (\sum\limits_{\tau' \in \F_{\tau,1}}\cH^1(I_{\tau'})) - \cH^1(I_\tau).
\end{gather} 
We see this as follows.  
By compactness of $\tau$ we have a point $\textrm{P}_1\in \bigcup\limits_{\tau' \in \F_{\tau,1}}I_{\tau'}$ such that  
\begin{gather*}
d_\tau=\dist(\textrm{P}_1,I_{\tau}).
\end{gather*}  
From the fact that the union is a union of line segments and $I_\tau$ is a line segment we have $\textrm{P}_1\in \tau$.
Let $\textrm{P}_2\in I_{\tau}$ satisfy $\dist(\textrm{P}_1,I_\tau)=\dist(\textrm{P}_1,\textrm{P}_2)$.
Set 
\begin{eqnarray*}
c_1&=&\dist(\textrm{P}_1,\tau(inital))\\
c_2&=&\dist(\textrm{P}_1,\tau(final))\\
a_1&=&\dist(\textrm{P}_2,\tau(initial))\\
a_2&=&\dist(\textrm{P}_2,\tau(final)).
\end{eqnarray*}
Hence (by say, the Cosine Theorem)
\begin{gather*}
d_\tau^2\leq c_1^2-a_1^2\\
d_\tau^2\leq c_2^2-a_2^2
\end{gather*}
and so 
\begin{eqnarray}
2d_\tau^2&\leq&c_1^2-a_1^2 + c_2^2-a_2^2\\
&=& (c_1-a_1)(c_1 +a_1) + (c_2-a_2)(c_2 +a_2)\\
&\leq& (c_1-a_1)2\diam(\tau) + (c_2-a_2)2\diam(\tau)\\ 
&=&2\diam(\tau)(c_1+c_2 - (a_1 + a_2)). 
\end{eqnarray}
Finally, by the triangle inequality, 
\begin{gather*}
c_1 + c_2 \leq \sum\limits_{\tau' \in \F_{\tau,1}}\cH^1(I_{\tau'})
\end{gather*}
and we also have 
\begin{gather*}
a_1 +a_2 = \cH^1(I_\tau)
\end{gather*}
which gives
\begin{gather*}
d_\tau^2 
\leq \diam(\tau)((\sum\limits_{\tau' \in \F_{\tau,1}}\cH^1(I_{\tau'})) - \cH^1(I_\tau))
\end{gather*}
which gives \eqref{temp2} as desired.

By summing \eqref{temp2} over $\F_n$ we have
\begin{gather*}
\sum\limits_{\tau \in \F_n}\frac{d_\tau^2}{\diam(\tau)} \lesssim
\sum\limits_{\tau \in \F_{n+1}}\cH^1(I_\tau)  -\sum\limits_{\tau \in \F_n}\cH^1(I_\tau).
\end{gather*}
Summing over all $n$ we get:  
\begin{gather*}
\sum\limits_{\tau \in \F}\frac{d_\tau^2}{\diam(\tau)} \lesssim
\supl_n(\sum\limits_{\tau \in \F_n}\cH^1(I_\tau)) \leq 
\length (\bigcup\limits_{\F_0} \tau).
\end{gather*}
We can now compute in an $\ell_2$ fashion:\\
\begin{eqnarray*}
(\sum\limits_{\tau \in \F} \bt(\tau)^2\diam(\tau))^\half &\leq&
  (\sum\limits_{\tau \in \F} \frac{(\sum\limits_{k=0}^\infty d_{\tau_k})^2}{\diam(\tau)})^\half\\
  & \leq&
  \sum\limits_{k=0}^\infty(\sum\limits_{\tau \in \F} \frac{d_{\tau_k}^2}{\diam(\tau)})^\half\\
  & \leq&
  \sum\limits_{k=0}^\infty2^{-J\frac{k}{2}}(\sum\limits_{\tau \in \F} \frac{d_{\tau_k}^2}{2^{-Jk}\diam(\tau)})^\half \\
&\lesssim&
  \sum\limits_{k=0}^\infty2^{-J\frac{k}{2}}(\sum\limits_{\tau \in \F} \frac{d_{\tau_k}^2}{ \diam(\tau_k)})^\half\\ &\lesssim&
  \sum\limits_{k=0}^\infty2^{-J\frac{k}{2}}(\length (\bigcup\limits_{\F_0} \tau))^\half\\
  & \lesssim&
  \length (\bigcup\limits_{\F_0} \tau)^\half
\end{eqnarray*}
where the penultimate inequality follows from the fact that $\tau \neq \tau' \then \tau_k \neq \tau'_k$
unless $\tau$ is one of the $\log(4A)$ immediate consecutive forefathers of $\tau'$ or vise-verse.
(More careful notation would eliminate this need for a factor of $\log(4A)$.) 
\end{proof}
\begin{rem}
If one follows the computation one gets
\begin{gather*}
\sum\limits_{\tau \in \F} \bt(\tau)^2\diam(\tau) \leq CA^\half \log(A) \length(\bigcup\limits_{\F_0} \tau),
\end{gather*}
where $C$ is a universal constant, independent of $A$.
\end{rem}

We now turn to the construction of filtrations.
As before, when we discuss the intersection or union of arcs, we do this in the parameter space (i.e. in $\T$\ ).  When we discuss the diameter of an arc, we do so in the image space (i.e. $\Gamma$\ ).

The idea for the proof of the lemma comes from now classical constructions of Dyadic Cubes on Homogeneous Spaces (see e.g. \cite{Ch}, and \cite{Da} page 93 for a simple version).  
One should note that condition (2) replaces a doubling condition.
\begin{lemma}\label{construct_filtration}
There is a universal constant $J>0$ ($J=10$ suffices) such that,  given a collection of arcs $\F^0=\bigcup\limits_{i=0}^\infty\F^0_i$  with the following properties:
\begin{verse}
(1) $\tau \in \F^0_{n} \then 
\sharp \{ \tau' \in \F^0_n: \tau' \cap \tau \neq \emptyset \} \leq C$\\
(2) $\tau \in \F^0_n \then 2^{-n} \leq \diam(\tau) \leq A2^{-n+1}$\ ,\\
\end{verse}
then we have (we construct) $2CJ$ families of arcs, each of which will be a filtration (see requirements of previous lemma).
Furthermore, we will have that for any $\tau' \in \F^0_{n'}$ there exists $\tau\in\F_{n}$ with $n'\sim Jn$ for one of the filtrations we construct. 
This $\tau$ will satisfy:
$\tau' \subset \tau$ and $\diam(\tau)<2\diam(\tau')$ 
(and hence $\bt(\tau)\geq c_0 \bt(\tau')$).
\end{lemma}
\begin{proof}
We will now construct $\leq (2CJ)$ filtrations: 
$\{\F^j\}_{j=1}^{k}, k \leq(2CJ)$.  
A single filtration will be denoted by $\F=\{\F_i\}_0^\infty$ (omitting the 
superscript) and will have the  properties required for the previous lemma.\\
We will use an order on the arcs given by the flow along a universally 
chosen parameterization of $\Gamma$.

First, divide each $\F^0_n$ into $C$ collections such that at every level $n$ each collection is composed of disjoint arcs.  Then divide each of these into 2 collections, such that within each collection any two 
arcs at  the same level $n$ will be separated by an arc of diameter at least $2^{-n}$.

Select a single collection from each level $n$ and call it $\F^1_n$.
Do NOT confuse this superscript with the enumeration of the different final filtrations; this is merely 
a step in the construction of a single filtration, $\F$.\\
Now, `dilute' each $\{\F^1_n\}_{n=0}^\infty$ by skipping $J$ generations at a time, 
multiplying the number of collections by $J$.  
Call a single collection $\F^2=\{F^2_n\}_{n=0}^\infty$ (renumbering $Jn \to n$).
 
We now want to turn $\{\F^2_n\}_{n=0}^\infty$ to a nested family.  
Consider $\tau^2 \in \F^2_n$.  
Set
\begin{eqnarray*} 
\tau_0&=&\tau^2\\
\tau_{k+1}&=&\tau_k \cup (\bigcup\limits_{\tau' \in \F^2\atop \tau_k \cap \tau' \neq \emptyset}\tau')\\ 
\tau^3&=&\lim\limits_k\tau_k.
\end{eqnarray*}
Denote   by $\F^3$  the family given by $\tau^2\to \tau^3$.

Note that $\bt$  (possibly) decreases by only a small factor (dependent on $J$) by $\tau^2\to \tau^3$, as 
the diameter increases by at most a factor of $1+4\cdot 2^{-J}$ (see Lemma 3.16 in \cite{RS-metric} for a simple proof of this by induction).
We have that $\F^3$ is almost the family which we desire (a filtration) -
requirement (4) and the existence part of (1) (see previous lemma) are not 
yet satisfied.

Suppose $n=0$.
Consider 
\begin{gather*}
R=\Gamma \setminus \bigcup\limits_{\tau \in \F^3_{n}}\tau = \cup R_j
\end{gather*} 
where $\{R_j\}$ are 
connected components, ordered by $\gamma$.  
Consider an $R_j$.  
Note that $\diam(R_j) > 2^{-n}$.
If \\
$\diam(R_j) > A2^{-n}$ 
then chop it up into a finite number of connected 
parts with $\diam \leq A2^{-n}$.
Rename them to be $\{R_j\}$.  
Go over them in order.  If an element $R_j$ has $\diam(R_j) < 2^{-n}$ then 
join it to the following element in $\F^3_n \cup\{R_j\}$.     
Now perturb each new $\partial R_j$ so that  elements of $\F^3_{n+1}$ have a father (unique).   
We call the collection of these new sets $\F^4_n$. (We remind the reader this was the all for $n=0$.)
This gives the requirements for the previous lemma for $n=0$.

Suppose we have the requirements for $n$ and we want to get them for $n+1$.
Consider $R=\Gamma \setminus 
(\bigcup\limits_{\tau \in \F^3_{n+1}}\tau  \cup \bigcup\limits_{\tau \in \F^4_{n}}\partial\tau)$. 
As before, we may write $R= \cup R_j$, where the $R_j$ are connected components.  Also, as before, we may subdivide  $R_j$ arcs to get arcs of diameter at most $A2^{-n-1}$.  We then rejoin them if necessary to  adjoining  $R_j$ arcs to make sure they are of diameter  at least $2^{-n-1}$. (Note that an $R_j$ arc must have been of diameter at least $2^{-n-1}$ before   being subdivided.)
By perturbing
each new $\partial R_j$ we make sure  $\F^3_{n+2}$ have a father (unique). 
We call the collection of these new sets $\F^4_{n+1}$.

We get that $\F^4$ is  the desired filtration $\F$.  
Clearly by making different initial choices we get a total of 
at most $2CJ$ filtrations.
 \end{proof}

\begin{lemma}
We have \eqref{beta_over_A_1} for  $j\in\{0,1,2\}$.
\end{lemma}
\begin{proof}
Fix $j\in\{0,1,2\}$.
We can fix $\epsilon>0$ below, independent of all other constants ($\epsilon$ will only serve us for the purpose of this lemma).

Notice that $\bt(\tau)$ is continuous in the endpoints of $\tau$.  
Also notice that if $n=-\log(\frac{\diam(2^jQ)}{A})$ and 
$\bt(\tau) \leq \epsilon$ for $\tau \in \Lambda(Q)$ then 
$\sharp(\tau \cap X_n) < 2^{j+1}A+1$.

If $\sharp(\gamma_{2^jQ} \cap X_n) \leq 2^{j+1}A +1$ set 
$\tau_Q = \gamma_{2^jQ}$.   
Otherwise, set $\tau_Q$ to be a  sub-arc of $\gamma_{2^jQ}$ 
such that it contains $center(Q)$ and
$\sharp(\tau_Q \cap X_n) = 2^{j+1}A +1$.
  
Using the definition of $\G_1$, we get that in both of the above cases  
\begin{gather*}
\bt(\tau_Q) \gtrsim \epsilon_2\beta(Q),\\
\sharp(\tau_Q \cap X_n) \leq 2^{j+1}A +1, \\
\diam(\tau_Q) \sim \diam(Q) \text{ (with constant $A$), and}\\
Q \to \tau_Q \text{ is at most } (2^{j+1}A +1):1.
\end{gather*}
Furthermore, 
\begin{gather}\label{041805}
\forall Q \in \G_1: 
\sharp\{Q' \in \G_1: \diam(Q) = \diam(Q'); \tau_{Q'} \cap \tau_Q \neq \emptyset\} 
\leq (3\cdot2^{j+1}A +3).
\end{gather}
To see this, we use the order on  $\{Q':\diam(Q')=\diam(Q)\}$,
given by the parameterization to the centers. 
Let $\tau_{Q'_1}$ be the largest, and $\tau_{Q'_2}$ the smallest  such elements. 
Then every other such $\tau_{Q'}$ must have the center of $Q'$ contained in  the union 
$\tau_{Q'_1} \cup \tau_{Q} \cup \tau_{Q'_2}$ which gives \eqref{041805}. 

For all of that we need $\epsilon$ fixed sufficiently small.
We now use the previous lemmas: first we use Lemma \ref{construct_filtration} and then Lemma \ref{L2_lemma}.\\
We have
\begin{gather*}
\sum\limits_{Q \in \G_1^j}\beta(Q)^2\diam(Q) \lesssim 
\sum\limits_{Q \in \G_1^j}\bt(\tau_Q)^2\diam(Q) \lesssim 
\sum\limits_{i=1}^C\sum\limits_{\tau \in \F^i}\bt(\tau)^2 \diam(\tau) \lesssim 
\length (\gamma)\lesssim \cH^1 (\Gamma).
\end{gather*}
\end{proof}
\begin{rem}
Note that if one follows the computation one gets that \eqref{beta_over_A_1} is satisfied  with constant
$\sim \frac{1}{\epsilon_2}^2 A^{\frac{5}{2}}\log(A)$. 
We will have (from subsection \ref{almost_straight_lines}) that  $\epsilon_2\sim {1\over A}$ and so we will get that \eqref{beta_over_A_1} is satisfied  with constant
$\sim  A^{\frac{9}{2}}\log(A)$. 
\end{rem}

\begin{lemma}
We have \eqref{beta_over_A_3}.
\end{lemma}
\begin{proof}
Consider $Q\in \G_3^0 \cap \G_3^1 \cap \G_3^2$.
From the definition of $\G_3^0$ we have the existence of 
$\xi^0_Q \in \Lambda_{Q}\smallsetminus S_Q$.
We consider two cases.
If $\xi^0_Q \subset \N_{\tenth \diam(Q)}(\gamma_Q)$,   
set $\tau_Q=\xi^0_Q$. 
If $\xi^0_Q \nsubseteq \N_{\tenth \diam(Q)}(\gamma_Q)$ then denote by 
$\xi^1_Q$ the extension of $\xi^0_Q$ to an element of $\Lambda_{2Q}$.
Set $\tau_Q=\xi^1_Q$.
In both cases we get 
\begin{gather*}
\diam(\tau_Q)\geq \half \diam(Q),
\end{gather*}
and by the definition of $\G_3^0$ and $\G_3^1$ 
(using the inequality $\beta(2Q) \geq \half \beta(Q)$) we have
\begin{gather*}
\bt(\tau_Q) \gtrsim \epsilon_2\beta(Q).
\end{gather*}
(This follows from the fact we may reduce $\epsilon_1$ in the definition of $\G_i^j$.)
We also have 
\begin{gather*}
\sharp \{Q' \in  \G_3^0 \cap \G_3^1 \cap \G_3^2: \diam(Q) = \diam(Q'); 
               \tau_Q \cap \tau_{Q'} \neq \emptyset\} 
\leq 8A  
\end{gather*}
which follows from  
$\gamma_{4Q} \supset \gamma_Q$ and $\beta_{S_{4Q}}(Q)$ being small.
Now we use  Lemma \ref{construct_filtration} and then Lemma \ref{L2_lemma} as we did 
in the proof of the previous lemma.

\end{proof}
\begin{rem}
Note that if one follows the computation one gets that \eqref{beta_over_A_3} is satisfied  with constant controlled by that of  equation \eqref{beta_over_A_1}.
\end{rem}
\begin{rem}[A remark concerning Lemma \ref{L2_lemma} and Lemma \ref{construct_filtration}.]
Consider a filtration $\F$.
The same proof as that of Lemma \ref{L2_lemma} gives us
\begin{gather}\label{BMO-statement}
\sum\limits_{\{\tau' \in \F; \tau' \subset \tau\}} \bt(\tau')^2\diam(\tau') 
\lesssim \length(\tau)
\end{gather}  
for all $\tau \in \F$.
Hence we can set 
\begin{gather*}
w_\tau(x)=\frac{\bt(\tau)^2\diam(\tau)}{\length(\tau)}, \forall x \in \tau
\end{gather*} 
and get from \eqref{BMO-statement}
\begin{verse}
(1) $\int_\tau w_{\tau}d\length= \bt(\tau)^2\diam(\tau)$\\  
(2) $w(x):= \sum\limits_{\tau \in \F} w_\tau(x)  \in BMO_{\F}$\\
(3) $supp (w_\tau) \subset \tau$
\end{verse}
where 
\begin{gather*}
f\in BMO_{\F} \iff 
\supl_{\tau\in\F} {1\over \length(\tau)}\int_\tau | f- {1\over \length(\tau)}\int _\tau f d\length | d\length < \infty. 
\end{gather*}
This is just a formal way of writing things, which is a useful reminder that this part of the non-Ahlfors-regular theory is close to the Ahlfors-regular  case. 
Another such reminder is Lemma \ref{construct_filtration}.  
One way of thinking about this is that even though $\Gamma$ is not Ahlfors-regular, we do have that $\T$ is Ahlfors-regular and so we may use standard ideas from the world of Ahlfor-regular theory (or homogeneous space theory).
\end{rem}
%

%
%
\subsection{Almost Flat Arcs}\label{almost_straight_lines}
In this subsection we prove \eqref{beta_over_A_2}.   
This is subsection is probably the hardest part of the paper, and so 
throughout this subsection,  the reader is urged to consider the example of $\Gamma$ being a 
finite union of straight line segments, ignoring any problems that may arise at the end-points of these segments .  
The proofs simplify somewhat if they are reduced to just this example, 
however almost all of the ideas will remain!

\begin{figure}[h]
\begin{center}
$\begin{array}{c}
\scalebox{0.3}{\includegraphics*[0in,0in][4in,4in]{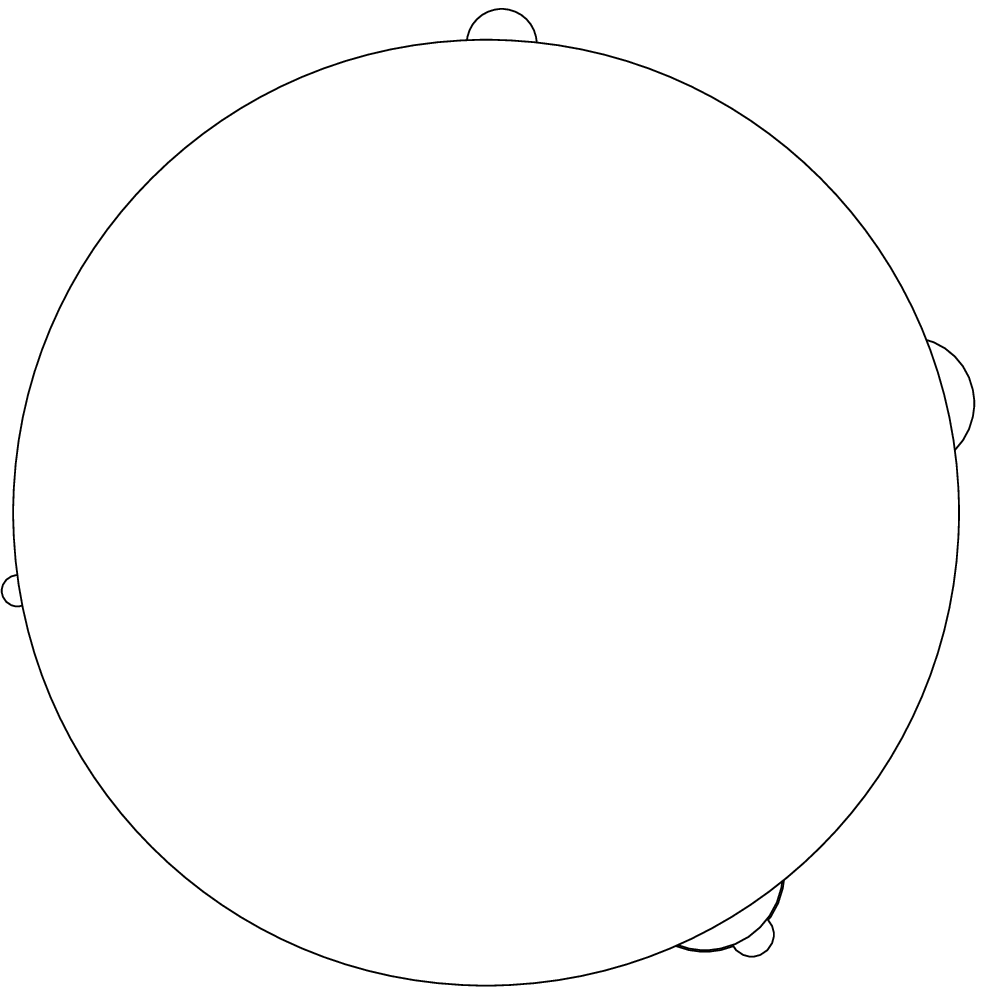}}
\end{array}$
\end{center}
\caption{Example of $U_Q$}
\label{U_Q-example}
\end{figure}

\begin{figure}[p]
\begin{center}
$\begin{array}{c@{\hspace{-.5in}}c@{\hspace{-.0in}}c}
\scalebox{0.1}{\includegraphics*[0in,0in][8in,8in]{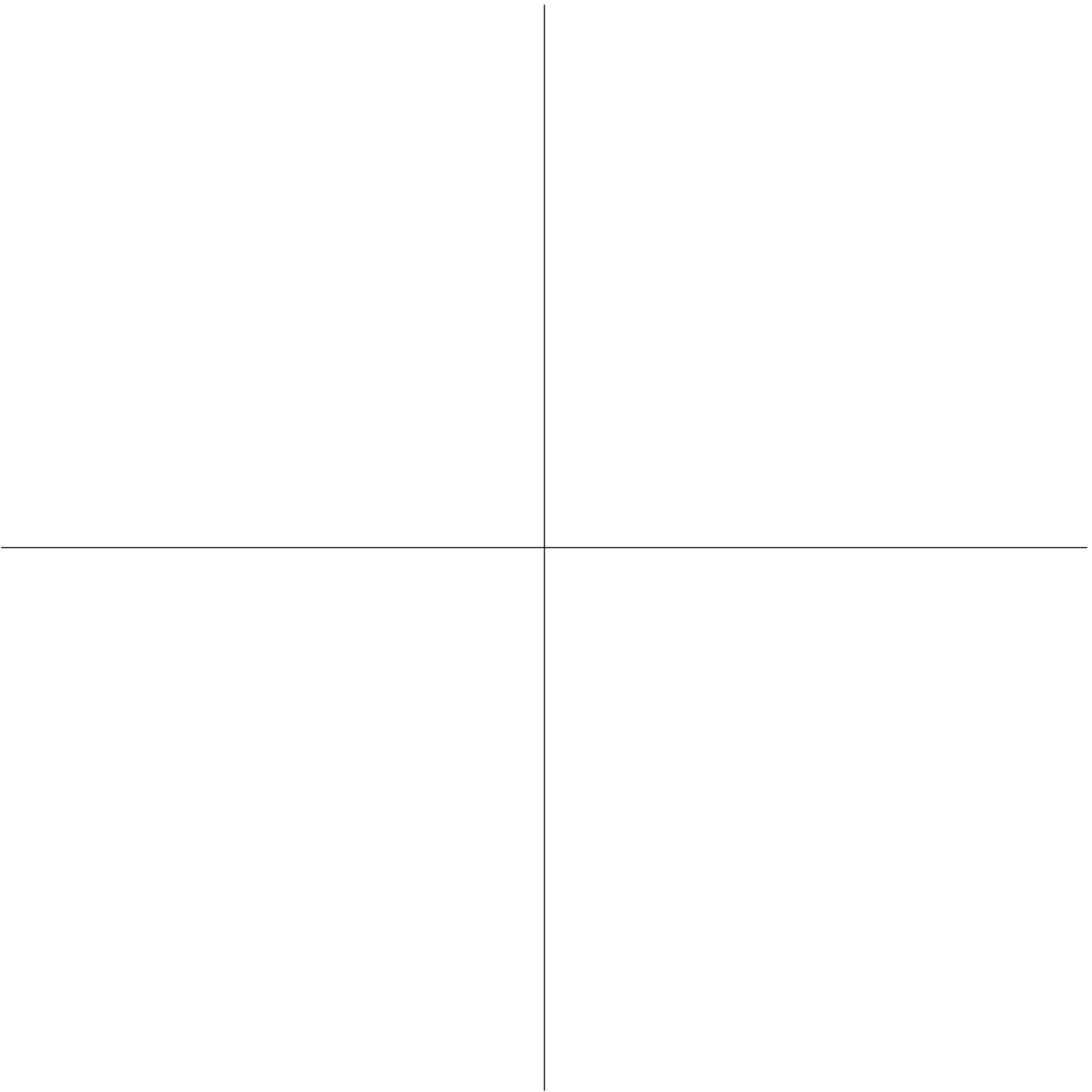}}&
\scalebox{0.1}{\includegraphics*[0in,0in][8in,8in]{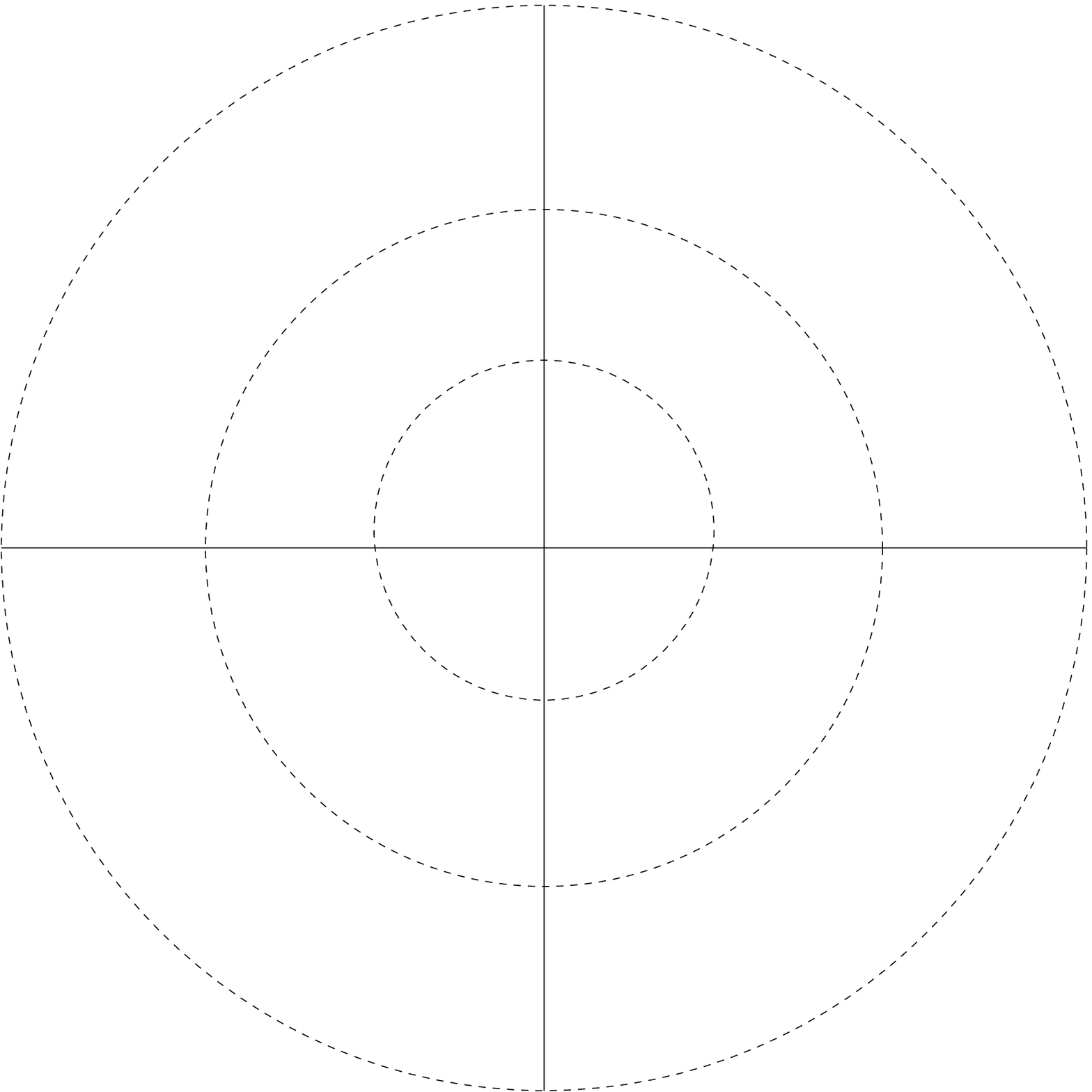}}&\\
\mbox{2 crossing segments.}&\mbox{\ \ \ \ \ \ \ \ \ \ \ Segments with $U_{Q_1}\supset U_{Q_2}\supset U_{Q_3}$ (circles).}&\\
&\mbox{We have $M=1$.}&\\
&\mbox{We assume these circles are the only $U_Q$'s.}&\\[0.4cm]
\scalebox{0.1}{\includegraphics*[0in,0in][8in,8in]{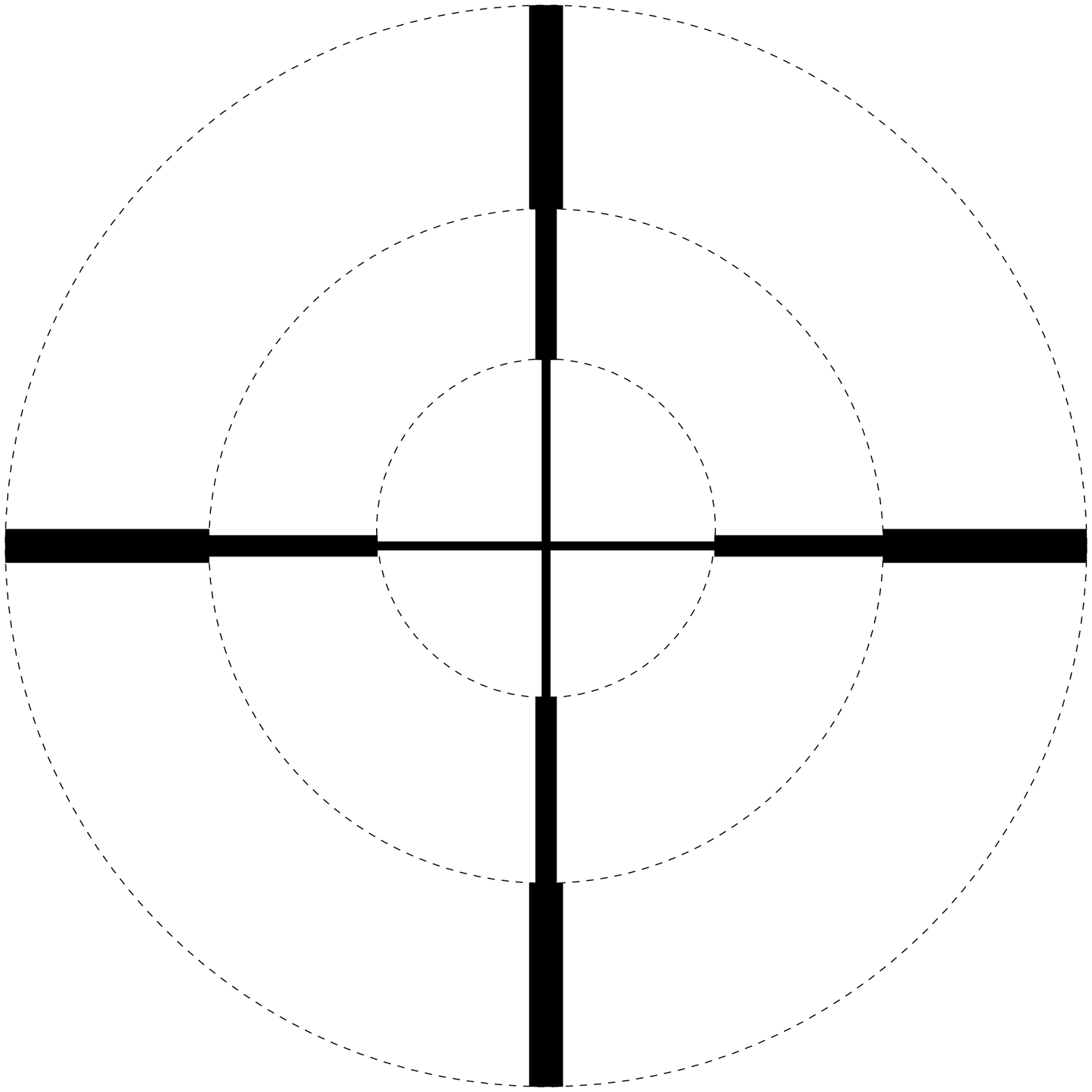}}&
   \hspace{-1.2in}\scalebox{0.1}{\includegraphics*[0in,0in][8in,8in]{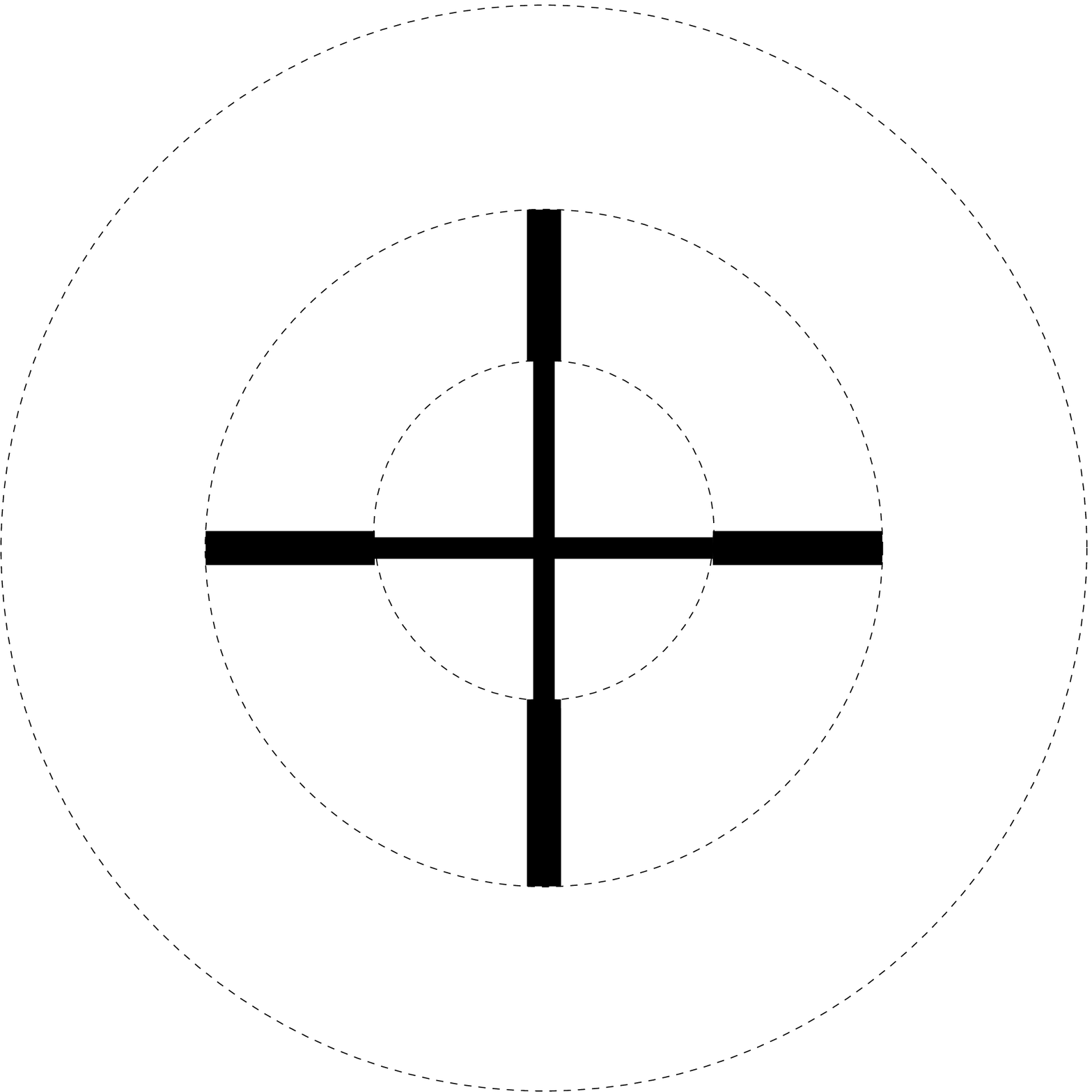}}&
   \hspace{-1.2in}   \scalebox{0.1}{\includegraphics*[0in,0in][8in,8in]{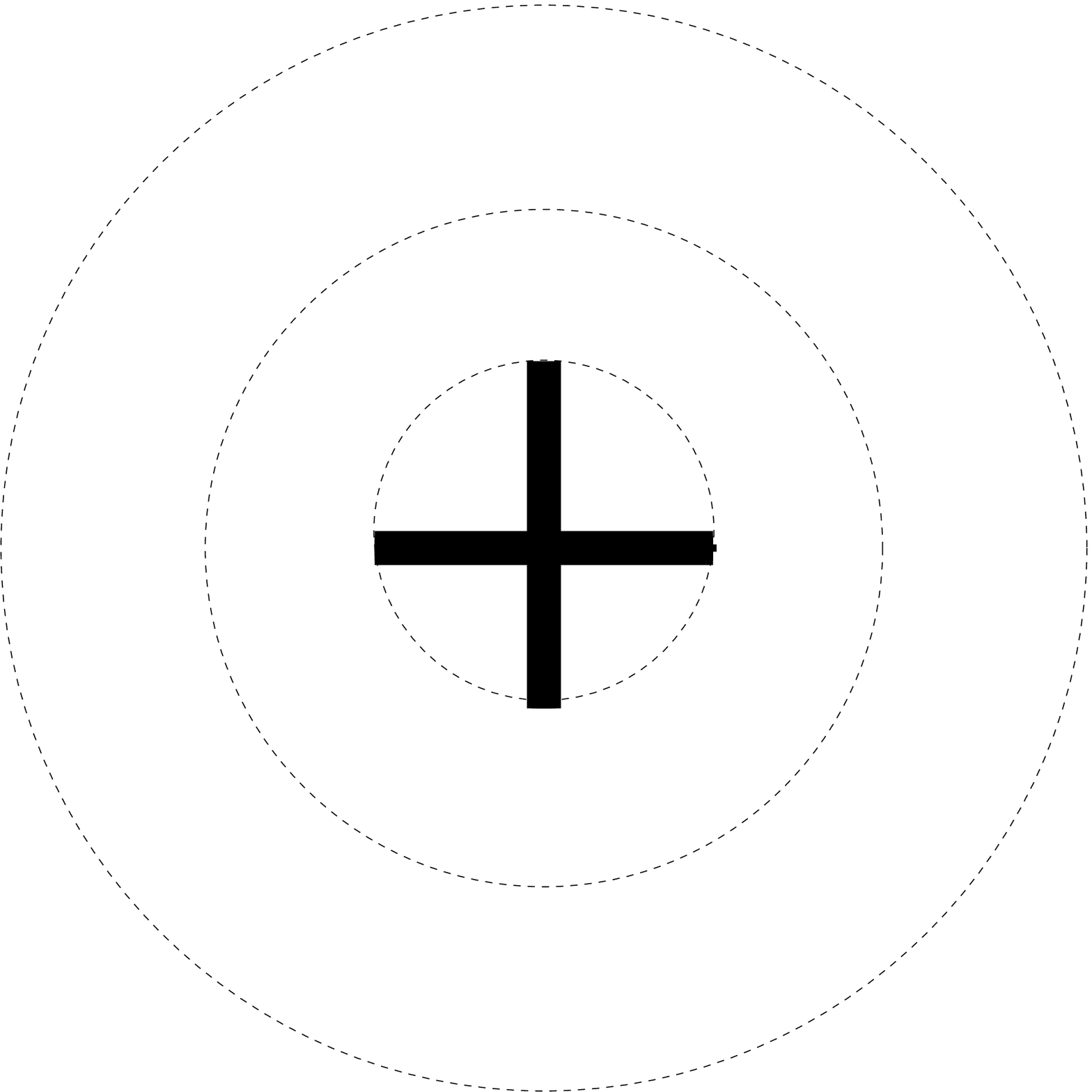}}\\
\mbox{$w_{Q_1}$}&\hspace{-1.2in}\mbox{$w_{Q_2}$}&\hspace{-1.2in}\mbox{$w_{Q_3}$}\\[0.4cm]
\scalebox{0.18}{\includegraphics*[0in,0in][8in,8in]{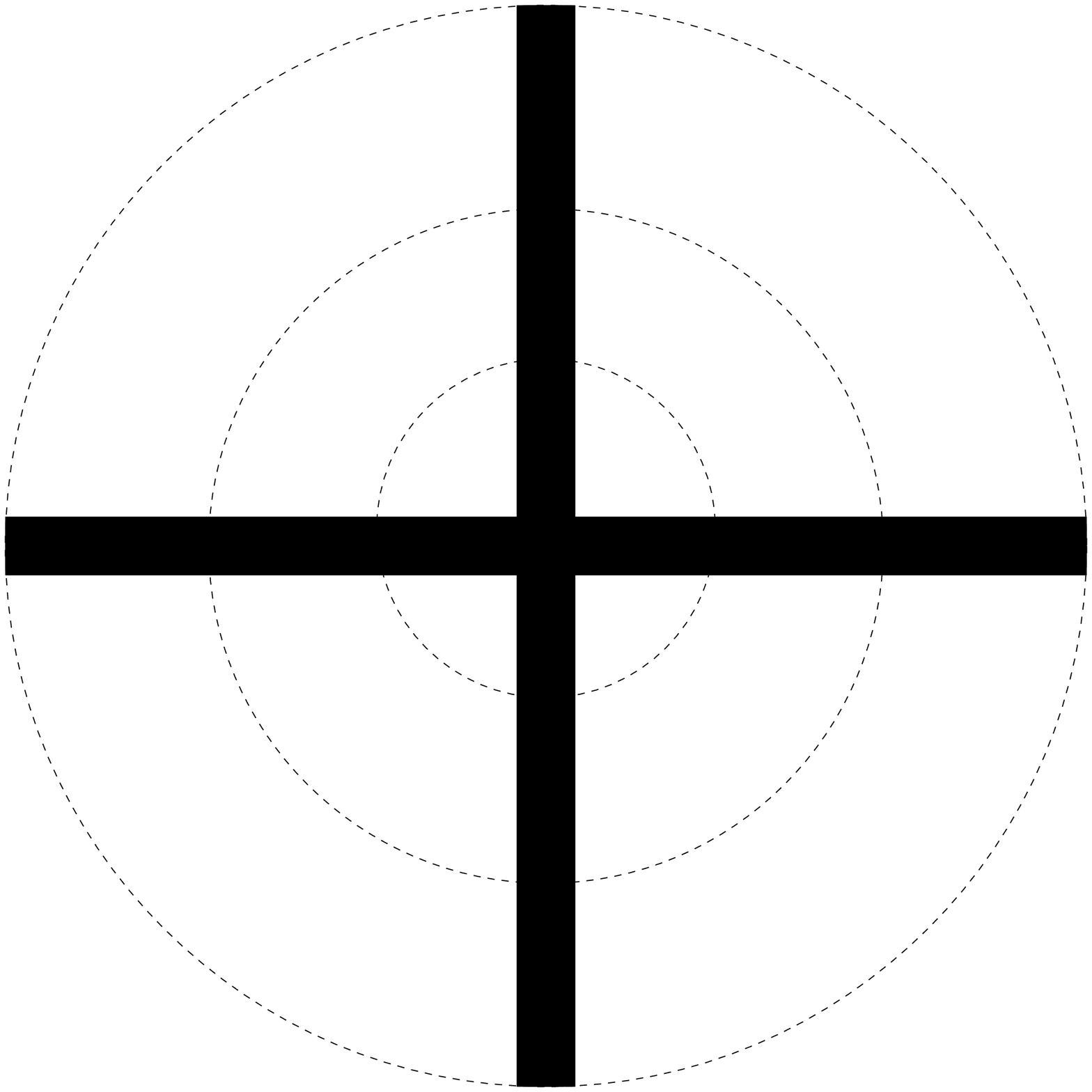}}&&\\
\mbox{The total weight = $w_{Q_1}+w_{Q_2}+w_{Q_3}$.}&&\\
\end{array}$
\end{center}
\caption{An  example of weight distribution arising from the martingale. 
Thickness of lines indicates value of $W_{Q_j}$.}
\label{weight}
\end{figure}

In \cite{Ok} Okikiolu proved a corresponding result by allotting for each cube $Q$ a segment $seg_Q$ whose length controlled $\beta(Q)\diam(Q)$ (Okikiolu  used the dyadic grid  as her multiresolution family). 
We follow in the same spirit by allotting a density (we use the word 'weight') for every ball.  
Let us give a vague idea of our plan:\\  
For each $Q$ we will define a density  (weight) function $w_Q$. (See Figure \ref{weight} for an example) 
We will have several families of balls (the number of such families is bounded by some universal constant).
Every ball $Q$ will have a \textit{core} $U_Q$ (see Figure \ref{U_Q-example}).
Within each family,  these \textit{cores} will have nice nesting properties between different balls.
We will get a constant $q<1$ such that if $U_{Q_1} \supset U_{Q_2}$ then 
$w_{Q_1}(x)\leq q \cdot w_{Q_2}(x)$.  
Hence, within each family, we will get that the sum of the densities at a given point is a geometric sum and hence bounded by a constant.  
To be slightly more accurate (in this vague setting), the above only happens for almost every point ($d\length$), which is enough.
Furthermore, $\int_Q w_Qd\length$ is enough to control the Jones beta number of the ball $Q$ (scaled correctly).\\

We start with a preliminary lemma.
We build for each ball $Q$ a \textit{core} $U_Q\subset Q$ such that these \textit{cores} are divided into $J$ families ($J$ being a sufficiently large universal constant).  Within each family they will have a nice 'nesting' structure (see property (4)). 

To see the origins of the idea for the statement  of Lemma \ref{building_U},  see  \cite{Da} page 93, or \cite{Ch}.
\begin{lemma}\label{building_U}
Given $c\leq \frac{1}{4A}$ and $J \geq 10$, there exist  $J$   families of connected sets in $H$
such that (denoting a single family by $\{U^{c,k}_n\}_{n=0,k=0}^{n=\infty,k=k_n}$):
\begin{verse}
       (1)  For every $x\in X_n$ there exists a unique $k$ such that $cQ\subset
       U^{c,k}_n$ for some family, where $radius(Q)=A2^{-n}$.\\
       (2)  $cA2^{-n} \leq \diam (U^{c,k}_n) \leq (1+4\cdot2^{-J+1})cA2^{-n}\,.$\\
       (3)  If $k \neq k'$ then $U^{c,k}_n \cap U^{c,k'}_n = \emptyset$  
as long as they are in 
the same family. In that case we also have $\dist( U^{c,k}_n,U^{c,k'}_n) \geq 2^{-n-1}$. \\
       (4)  If $U^{c,k}_m \cap U^{c,k'}_n \neq \emptyset$, they are in the same family and $m>n$, then 
            $U^{c,k}_m \subset U^{c,k'}_n$.
\end{verse}
\end{lemma}
\begin{proof}
Let  $Q_n^k=\ball(x^k,A2^{-n})$ where $x^k \in X_n$ for the 
proof of this lemma.  Then if $k\neq k'$ then $\dist (cQ_n^k, cQ_n^{k'}) > 2^{-n}$.
Set 
\begin{eqnarray*}
U^{c,k}_{m,0}&:=& cQ_m^k\\
U^{c,k}_{m,i+1}&:=& U^{c,k}_{m,i} \cup \bigcup\limits_{cQ_{m+i'J}^{k'}\cap U^{c,k}_{m,i} \neq \emptyset} cQ_{m+i'J}^{k'}\\
U^{c,k}_{m}&:=&\lim U^{c,k}_{m,i}
\end{eqnarray*}

The $j^{th}$ family is
\begin{gather*}
\{U^{c,k}_{m}: m\in j + J\mathbb{N}\}.
\end{gather*}
See Lemma 3.16 in \cite{RS-metric} for proof of property (2). The other properties easily follow from the definitions.
\end{proof}

\begin{rem}
When changing $c$ we change only the parental relationship (i.e. the tree structure) within each family, and when we change $J$ we mix between families.  
\end{rem}
Set
\begin{eqnarray*}
c_0&:=&\frac {1}{64A}\\
U_Q&:=&U^{c_0,k}_n\\
U^\ceight_Q&:=&U^{8c_0,k}_n\\
U^\cfour_Q&:=&U^{16c_0,k}_n
\end{eqnarray*}
where $Q=\ball(x^k,A2^{-n})$ and $x^k\in X_n$.

The purpose of the following is to prove  \eqref{beta_over_A_2} for $j\in\{0,1,2\}$. 
We will have
\begin{gather}\label{041805-2}
\sum\limits_{Q \in \G_2^j}\beta(Q)^2\diam(Q) \lesssim  
\sum\limits_{Q \in \G_2^j}\beta_{S_{2^jQ}}(2^jQ)^2\diam(2^jQ) \leq  
C \length (\gamma).
\end{gather}  
All but the last inequality are obvious.\\ 
We will show
\begin{proposition}\label{sum_straight_prop}
\begin{gather}\label{sum_straight_eq}
\sum\limits_{Q \in \G_2}\beta_{S_{Q}}(Q)\diam(Q) \leq C \cH^1 (\Gamma).
\end{gather}
\end{proposition}
\noindent 
Since  $A$ is arbitrary, this will give us the last inequality in \eqref{041805-2} and hence \eqref{beta_over_A_2} for $j\in\{0,1,2\}$.  

\begin{figure}[p]
\begin{center}
$\begin{array}{c}
\scalebox{0.5}{\includegraphics*[1in,1in][8in,8in]{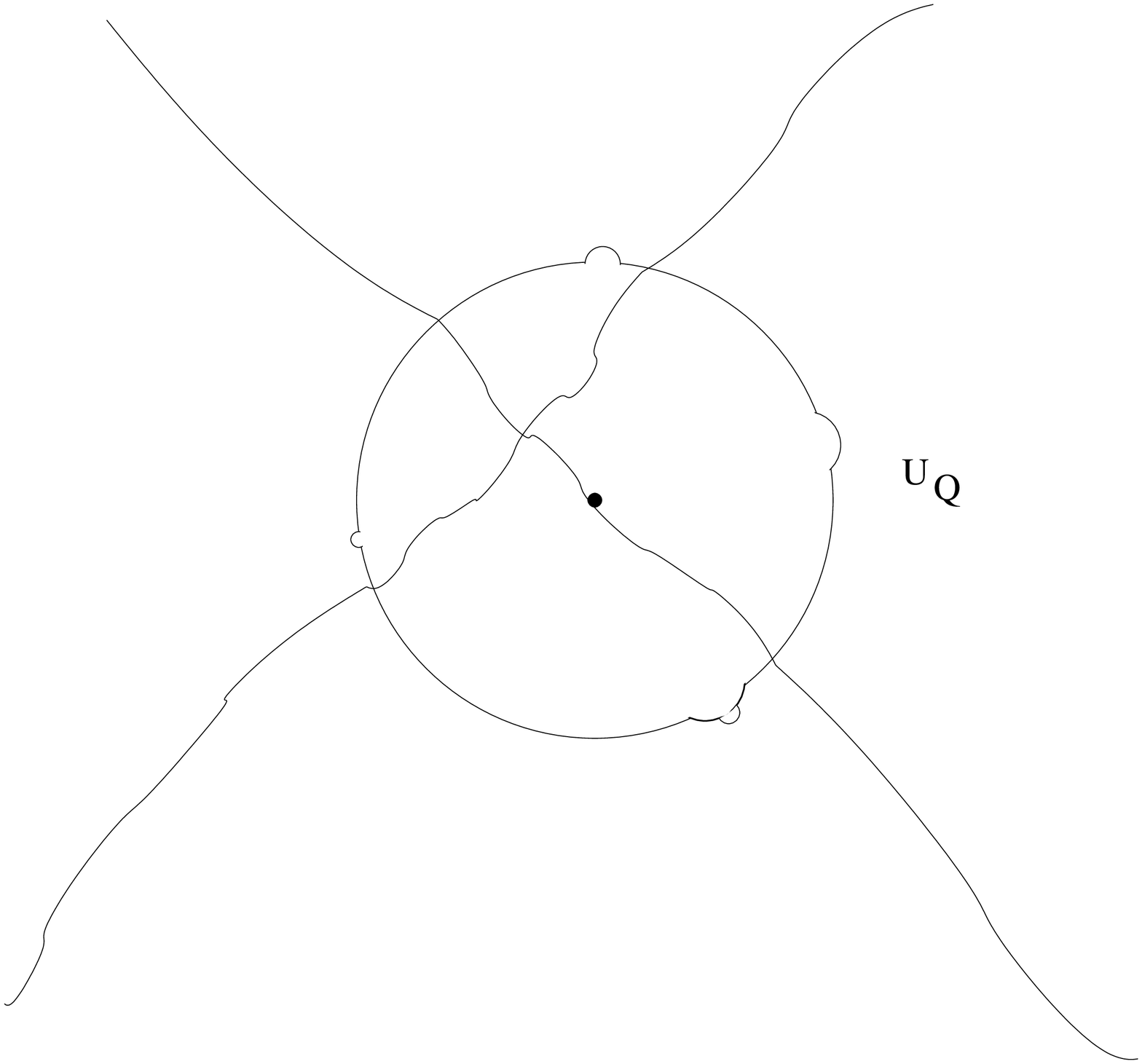}}\\
\vspace{.3in}\\
\scalebox{0.5}{\includegraphics*[0in,0in][8in,8in]{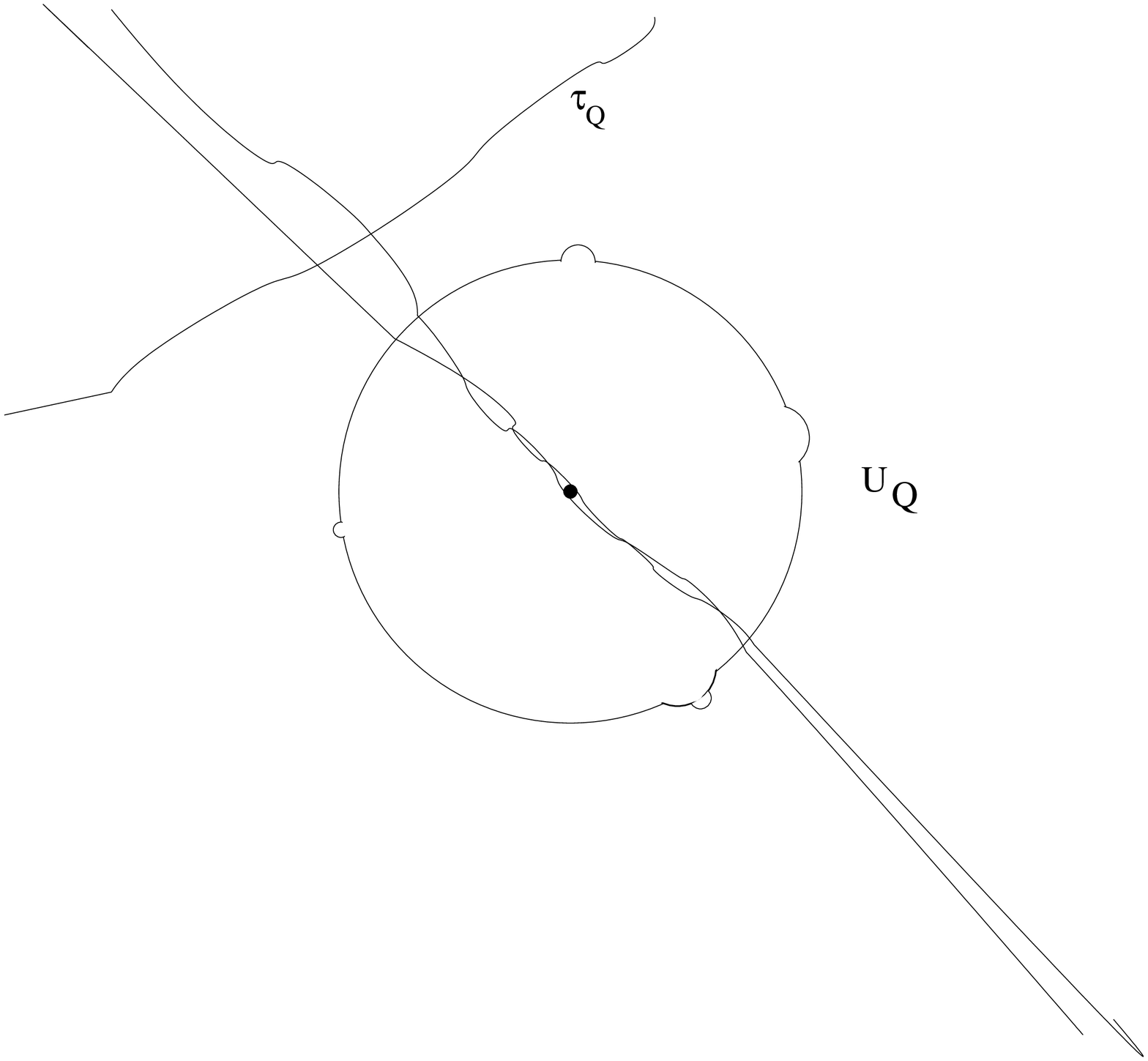}}\\
\end{array}$
\end{center}
\caption{Example of a $\Delta_1$ element  (top) and a $\Delta_{2}$ element (bottom)}
\label{U-examples}
\end{figure}

\begin{rem}
One should note the lack of the power $2$ in equation \eqref{sum_straight_eq}.
\end{rem}
Set (for a constant $C_U$ which can be fixed at the end of the proof) 
\begin{eqnarray*}
\Delta_1&=&\{Q \in \G_2: 
C_U\beta_{S_Q}(U^\ceight_Q) > \beta_{S_Q}(Q)\}\\
\Delta_2&=& \G_2 \smallsetminus \Delta_1
\end{eqnarray*}
See Figure \ref{U-examples} for  examples.
\begin{rem}\label{09042005}
We have (taking $C_U$ large enough) for all $Q \in \Delta_2$\\
(a.)\quad  $\beta_{S_Q}(4c_0Q) < \epsilon_0$. 
	(If we want $\epsilon_0$ small independent of $A$ then  we need $C_U\sim A$.) \\ 
(b.)\quad  the existence of 
	$\tau_Q \in S_Q$ such that $\tau_Q \cap U^\csixtyfour_Q = \emptyset$.\\
(c.)\quad (b) implies that we have  that  $\beta_{S_Q}(Q) \gtrsim 1$ for $Q\in \Delta_2$.
\end{rem}

\noindent
Set
\begin{eqnarray}
\Delta_{2.2}&=&\{Q\in \Delta_2:
	(\Lambda(Q)\smallsetminus S_Q) \cap U^\ceight_Q \neq \emptyset \}\\
\Delta_{2.1}&=&\Delta_2 \smallsetminus \Delta_{2.2}\,.
\end{eqnarray}

\noindent
We first show control over 
\begin{gather*}
\sum\limits_{Q \in \Delta_{2.2}}\beta_{S_{Q}}(Q)\diam(Q).
\end{gather*}
This is a slightly stronger version of what is done in section subsection \ref{curvy_arcs}  for the proof of \eqref{beta_over_A_1}.

\begin{lemma}\label{delta_2.2}
\begin{gather*}
\sum\limits_{Q \in \Delta_{2.2}}\beta_{S_{Q}}(Q)\diam(Q) \lesssim
\length (\gamma).
\end{gather*}
\end{lemma}
\begin{proof}
Let $\xi^0_Q\in \Lambda(Q)\smallsetminus S_Q$ such that 
\begin{gather*}
\xi^0_Q \cap U^\ceight_Q \neq \emptyset.
\end{gather*}
This implies 
\begin{gather*}
\diam (\xi^0_Q \cap U^\cfour_Q) > 8c_0 \diam(Q).
\end{gather*}
Suppose 
\begin{gather*}
\xi^0_Q \cap U^\ceight_{Q_i} \neq \emptyset
\end{gather*} 
for different balls 
$Q_i \in \Delta_{2.2}$ of the same diameter as $\diam(Q)$.
We have 
\begin{gather*}
\diam (\xi^0_Q \cap U^\cfour_{Q_i}) > 8c_0 \diam(Q)
\end{gather*} 
for all but maybe 2 (on $\partial \xi^0_Q$).
If we have an arc of diameter $1$ decomposed into $n$ disjoint sub-arcs, such that all but maybe two are of diameter $8c_0$ 
then $n>9A$ implies the arc has $\beta \geq \epsilon$.

Hence we can do as follows:  If $\xi^0_Q$ intersects $\leq 9A$ other $U^\ceight_{Q_i}$'s of 
the same scale, then set $\xi_Q=\xi^0_Q$.  Otherwise take a sub-arc which intersects exactly $9A$ 
other $U^\ceight_{Q_i}$'s of the same scale (and also intersects our $U^\ceight_Q$)
to be $\xi_Q$.
We have
\begin{gather*}
\bt(\xi_Q) \gtrsim \epsilon_2\beta(Q),
\end{gather*}
\begin{gather*}
\diam(\xi_Q) \sim \diam(Q) \text{ (with constant $A$)}. 
\end{gather*}
Furthermore, 
\begin{gather*}
\sharp\{Q' \in \Delta_{2.2}: \diam(Q) = \diam(Q'); \xi_{Q'} \cap \xi_Q \neq \emptyset\} 
\leq 9A.
\end{gather*}
We are now in position to use the lemma's of subsection \ref{curvy_arcs}.  
First we use Lemma \ref{construct_filtration} and then Lemma \ref{L2_lemma}.
This gives us
\begin{gather*}
\sum\limits_{Q \in \Delta_{2.2}}\beta_{S_{Q}}(Q)^2\diam(Q) \lesssim
\length (\gamma).
\end{gather*}
By Remark \ref{09042005} (c) this gives us the lemma.
\end{proof}
%
%
We now turn to deal with $\Delta_{2.1}$ and $\Delta_1$.

\subsubsection{Summing over $\Delta_{1}$}\label{delta-1-subsub}
Fix $M\in \mathbb{N}$ (we will sum over $M\geq 0$ in  corollary \ref{cor:delta_1-sum}).
Define $\Delta$ as 
\begin{gather*}
\Delta=\Delta(M):=\{Q\in \Delta_1: 2^{-M} \leq \beta(U^\ceight_Q) < 2^{-M+1} \}.
\end{gather*}
Take $K$ such that $1\leq K < MJ$.
Let $\Delta' \subset \Delta$,  be such that 
\begin{gather*}
\Delta'=\Delta'(M,K):=\{Q\in \Delta: \radius(Q)\in A2^{K+MJ\mathbb{N}} \}.
\end{gather*}
In other words, $\Delta'$ is obtained from  $\Delta$  by  {\it thinning} it, i.e. by taking every \{MJ\}-{th} element (starting at some offset $K$).

Consider $Q \in \Delta'$.  
Write 
\begin{gather}\label{20_09_05}
U^\cfour_Q\cap \Gamma=(\bigcup\limits_i U^\cfour_{Q^i}\cap \Gamma) \cup R_Q
\end{gather} where
$U^\cfour_{Q^i}$ is maximal in $U^\cfour_Q$, such that $Q^i \in \Delta'$ and
$R_Q=\Gamma \cap U^\cfour_Q \smallsetminus \bigcup\limits_i U^\cfour_{Q^i}$.
(By $U^\cfour_{Q^i}$ `maximal'  in $U^\cfour_Q$ we mean that there does not exist $Q'\in \Delta'$ such that 
$U^\cfour_{Q^i}\subset U^\cfour_{Q'} \subset U^\cfour_Q$.) 

\begin{lemma}\label{delta_1}
Suppose $\Delta'=\Delta'(M,K)$ is as above. 
Then 
\begin{gather*}
\sum\limits_{Q \in \Delta'}\beta(Q)\diam(Q) \leq C_4 2^{-M}\length (\gamma).
\end{gather*} 
\end{lemma}
\begin{proof}
We will construct weights that satisfy (i), (ii) and (iii):\\  
\indent(i) $\int_Q w_Qd\length> C_1 \diam(Q)$.\\
\indent(ii) for almost every $x_0\in \Gamma , \sum\limits_{Q \in \Delta'} w_Q(x) < C_2$.\\
\indent(iii) $\supp (w_Q) = U^\cfour_Q \cap \Gamma$.\\
This will suffice as then (using (i) and (ii))
\begin{gather*}
\sum\limits_{Q \in \Delta'}\beta(Q)\diam(Q) \lesssim 
\sum\limits_{Q \in \Delta'}2^{-M}\int_Q w_Qd\length  \leq 
2^{-M}\int_\Gamma \sum\limits_{Q \in \Delta'} w_Qd\length \lesssim
2^{-M}\int_\Gamma d\length \leq
2^{-M}\length (\gamma)\,,
\end{gather*} 
giving the lemma.
Consider now $Q \in \Delta'$.  
We construct $w_Q$ as a martingale.  
We will write $w_Q(V)$ for $\int_V w_Qd\length$ for any (measurable) set $V$.\\ 

\noindent
Set
\begin{gather*}
w_Q(U^\cfour_Q):=\diam(U^\cfour_Q ).
\end{gather*}
Given $w_Q(U^\cfour_{Q'})$, where
\begin{gather*}
U^\cfour_{Q'}\cap \Gamma=(\bigcup U^\cfour_{Q'^{j}}\cap \Gamma) \cup R_{Q'}
\end{gather*} 
as in equation \eqref{20_09_05},
then set
\begin{gather*}
w_Q(R_{Q'}):=\frac{w_Q(U^\cfour_{Q'})}{s'} \length(R_{Q'})
\end{gather*}
and 
\begin{gather*}
w_Q(U^\cfour_{Q'^{j}}):=\frac{w_Q(U^\cfour_{Q'})}{s'}\diam(U^\cfour_{Q'^{j}})
\end{gather*}
where 
\begin{gather*}
s':=\length(R_{Q'}) + \sum\limits_j \diam(U^\cfour_{Q'^{j}}).
\end{gather*}

Note that if we have an arc $\xi\subset\hat{\xi}\cap U^\cfour_{Q'}$ where $\hat{\xi}\in \Lambda(Q')$ then 
\begin{gather*}
\diam(\xi) \leq \length(R_{Q'}\cap \xi) + \sum\limits_{U^\cfour_{Q'^j}\cap \xi \neq \emptyset} \diam(U^\cfour_{Q'^j}).
\end{gather*}
Also note that 
\begin{gather*}
s'\leq (1 + 2^{-J+1})\length(\Gamma\cap U^\cfour_Q)< \infty
\end{gather*}
by considering 
$\gamma_{Q'^j}\cap U^\cfour_{Q'^j}$.
Now,\\
\noindent{\bf Step 1:}

There exists a universal constant $q<1$ such that 
\begin{gather*}
\frac{\diam(U^\cfour_{Q'})}{s'} \leq q.
\end{gather*}
To see this, 
notice   
we have an arc $\xi^0\subset S_{Q'}$ such that 
\begin{gather*}
\beta_{\xi^0 \cup \gamma_{Q'}}(U^\ceight_{Q'})>\half \beta(U^\ceight_{Q'})\,.
\end{gather*}
The above follows trivially in the case all arcs in $S_{Q'}$ are straight line segments.  The general case follows  from 
$\beta_{S_{Q'}}(U^\ceight_{Q'})\ge C_U^{-1}\beta_{S_{Q'}}(Q')\ge
 C_U^{-1}\epsilon_1\beta(Q')$,  the definition of $S_{Q'}$, and ensuring $\epsilon_2$ is sufficiently small with respect to $C_U^{-1}\epsilon_1$.
 
Now, let $\eta_{Q'}$, be a largest connected component of $\gamma_{Q'}\cap U^\cfour_{Q'}$.
Let $\xi$ be the largest connected component of $\xi^0\cap U^\cfour_{Q'}$.
By considering both $\xi$ and $\eta_Q$ we get 
\begin{eqnarray*}
\frac{\diam(U^\cfour_{Q'})}{s'}
&=&
\frac{\diam(U^\cfour_{Q'})}{\length(R_{Q'}) + \sum\limits_j \diam(U^\cfour_{Q'^{j}})}\\
&\leq&
{\diam(U^\cfour_{Q'}) \over
	\length(R_{Q'}) + 
	 \sum\limits_{U^\cfour_{Q'^{j}} \cap \eta_{Q'}\neq \emptyset} \diam(U^\cfour_{Q'^{j}}) +
	 \sum\limits_{U^\cfour_{Q'^{j}} \cap \xi \neq \emptyset \atop
	 			U^\cfour_{Q'^{j}} \cap \eta_{Q'}= \emptyset} 
	 	\diam(U^\cfour_{Q'^j})                }\\
&\leq& 
{(1 + 2^{-J+1})16c_0\diam(Q') \over
	\diam(\eta_{Q'}) + 
	\length(R_{Q'}\smallsetminus  \eta_{Q'}) + 
	\sum\limits_{U^\cfour_{Q'^{j}} \cap \xi \neq \emptyset \atop
		U^\cfour_{Q'^{j}} \cap \eta_{Q'}= \emptyset}  \diam(U^\cfour_{Q'^{j}})                }\\
&\leq& 
{ (1 + 2^{-J+1})16c_0\diam(Q') \over  
	16c_0\diam(Q') +  {1\over 10} c_0\diam(Q') }\\
&=& 
{ (1 + 2^{-J+1}) \over  
	1  +  {1\over 160}} \\
&\leq& q.
\end{eqnarray*}
where $q<1$ is a universal constant. (As we may impose $J>10$.)

\noindent{\bf Step 2:}

We now have 
\begin{eqnarray*}
\frac{w_Q(U^\cfour_{Q'^{j^*}})}{\diam(U^\cfour_{Q'^{j^*}})} 
&=&\frac{w_Q(U^\cfour_{Q'})}{s'}\\
&=&\frac{w_Q(U^\cfour_{Q'})}{\diam(U^\cfour_{Q'})} \frac{\diam(U^\cfour_{Q'})}{s'}\\
&\leq& 
 	q   \frac{w_Q(U^\cfour_{Q'})}{\diam(U^\cfour_{Q'})} 
\end{eqnarray*}
where $q<1$ is a universal constant.

\noindent{\bf Step 3:}

We observe that \textit{Step 2} gave us more.  Suppose now  that  $x\in U^\cfour_{Q_N} \subset ...\subset U^\cfour_{Q_1}$.
Using \textit{step 2} with 
\begin{gather*}
Q=Q_1;\qquad Q'=Q_n;\qquad Q'^{j*}=Q_{n+1}
\end{gather*}
we  get:
\begin{eqnarray*}
\frac{w_{Q_1}(U^\cfour_{Q_N})}{\diam (U^\cfour_{Q_N})} &\leq& 
  q\frac{w_{Q_1}(U^\cfour_{Q_{N-1}})}{\diam (U^\cfour_{Q_{N-1}})} \\
  &\leq&
  q^2\frac{w_{Q_1}(U^\cfour_{Q_{N-2}})}{\diam (U^\cfour_{Q_{N-2}})} \leq...\leq
  q^{N-1}\frac{w_{Q_1}(U^\cfour_{Q_{1}})}{\diam (U^\cfour_{Q_{1}})}=q^{N-1}.
\end{eqnarray*}
 Hence,  we have $w_{Q_1}(x) \leq 2q^{N-1}$, 
 and so we have (ii) as a sum of a geometric series.
\end{proof}

\begin{cor}\label{cor:delta_1-sum}
\begin{eqnarray*}
\sum\limits_{Q \in \Delta_1}\beta_{S_{Q}}(Q)\diam(Q)&\leq&
\sum\limits_{M\geq0}\sum\limits_{Q \in \Delta_1\atop2^{-M-1} < \beta_{S_{Q}}(Q) \leq 2^{-M}} \beta_{S_{Q}}(Q)\diam(Q) \\
&\leq& C_U J\sum\limits_{M\geq0}M2^{-M} \length (\gamma)\\
&\lesssim&
\length (\gamma)
\end{eqnarray*}
\end{cor}

%
%

\subsubsection{Summing over $\Delta_{2.1}$}
Consider a set $\Delta' \subset \Delta_{2.1}$, such that $\Delta'$ contains only balls from a single family as  constructed in Lemma \ref{building_U}.
Consider $Q \in \Delta'$.  
Write 
\begin{gather}\label{090605-1}
U_Q\cap \Gamma=(\bigcup\limits_i U_{Q^i}\cap \Gamma) \cup R_Q
\end{gather} where
$U_{Q^i}$ is maximal in $U_Q$, such that $Q^i \in \Delta'$ and
$R_Q=\Gamma \cap U_Q \smallsetminus \bigcup\limits_i U_{Q^i}$.
For a given $j$, We will denote the continuations of $\tau_{Q^{j}}$ and $\gamma_{Q^{j}}$ to arcs
in $\Lambda(Q)$ by 
$\widehat{\tau_{Q^{j}}}$ and $\widehat{\gamma_{Q^{j}}}$ respectively 
(we remind the reader of the $\tau$  assured by Remark \ref{09042005}.b).
\begin{rem}\label{r:all-is flat-and-lots-of-length}
A key observation we will use is that  
if $J$ (from Lemma \ref{building_U}) is large enough ($J\sim\log(A)$ suffices) we have
(for $\widehat{\tau_{Q^{j}}}$ as defined above)
\begin{gather*}
\widehat{\tau_{Q^{j}}}\in S_Q, 
\end{gather*}
since otherwise we would have had $Q\in \Delta_{2.2}$.
We also have 
$\widehat{\tau_{Q^{j}}} \cap 3c_0Q \neq \emptyset$.
Combining the two we get
$\widehat{\tau_{Q^{j}}} \cap c_0Q \subset \N_{\epsilon_0 4c_0\diam(Q)}(\gamma_{Q})$  and $\diam(\widehat{\tau_{Q^{j_0}}} \cap U_{Q}) > (1-8\epsilon_0) \diam(U_{Q}).$ 
($\epsilon_0$ is defined in  Remark \ref{09042005}.a)
Similarly , $\widehat{\gamma_{Q^{j}}}$ has the same properties.
\end{rem}
\begin{lemma}\label{delta_2.1}
Suppose $\Delta'$ is as above. 
Then 
\begin{gather*}
\sum\limits_{Q \in \Delta'}\beta(Q)\diam(Q) \leq C_4 \length (\gamma).
\end{gather*} 
\end{lemma}
\begin{proof}
We will construct weights that satisfy (i), (ii) and (iii):\\  
\indent(i) $\int_Q w_Qd\length> C_1 \beta(Q)\diam(Q)$\\
\indent(ii) for almost every $x\in \Gamma,  \sum\limits_{Q \in \Delta'} w_Q(x) < C_2$\\
\indent(iii) $\supp (w_Q) = U_Q \cap \Gamma$.\\
As in the proof of Lemma \ref{delta_1}, this is sufficient to give Lemma \ref{delta_2.1}.

Consider $Q \in \Delta'$.  
We construct $w_Q$ as a martingale.  
We do so in a similar manner to what is done in  Lemma \ref{delta_1}. We must be more careful here.

\noindent
Set
\begin{gather*}
w_Q(U_Q):=\diam(U_Q ).
\end{gather*}
Given $w_Q(U_{Q'})$, where
\begin{gather*}
U_{Q'}\cap \Gamma=(\bigcup U_{Q'^{j}}\cap \Gamma) \cup R_{Q'}
\end{gather*} 
as in equation \eqref{090605-1},
then
\begin{gather}\label{e:2_in_R_Q}
w_Q(R_{Q'}):=\frac{w_Q(U_{Q'})}{s'}  \length(R_{Q'}) \cdot 2
\end{gather}
and 
\begin{gather*}
w_Q(U_{Q'^{j}}):=\frac{w_Q(U_{Q'})}{s'} \diam(U_{Q'^{j}}) 
\end{gather*}
where 
\begin{gather*}
s':=\length(R_{Q'})\cdot 2 + \sum\limits_j  \diam(U_{Q'^{j}}).
\end{gather*}
Note that 
\begin{gather*}
s'\leq 2(1 + 2^{-J+1})\length(\Gamma\cap U_Q)< \infty
\end{gather*}
by considering 
$\gamma_{Q'^j}\cap U_{Q'^j}$.
Now,\\
\noindent{\bf Step 1:}

There exists a universal constant $q<1$ such that 
\begin{gather*}
\frac{\diam(U_{Q'})}{s'} \leq q.
\end{gather*}
To see this, let $\eta=\eta_{Q'}$ be a largest connected component of $\gamma_{Q'}\cap U_{Q'}$.  
We know 
\begin{gather*}
(1 + 4\cdot 2^{-J+1}) \diam (\eta_{Q'}) \geq \diam(U_{Q'}) \geq \diam (\eta_{Q'})\geq (1-8\epsilon_0)c_0\diam(Q).
\end{gather*}
(Recall that $\epsilon_0$ is defined in  Remark \ref{09042005}.a.)
Consider $I_\eta:=[\eta(initial),\eta(final)]$. 
Let $\pi:\Gamma\cap(1-4\epsilon_0)c_0Q' \to I_\eta$ be the radial (orthogonal) projection.  
If 
\begin{gather*}
U_{Q'^j} \subset (1-2^{-J+1}-4\epsilon_0)c_0Q'
\end{gather*} 
then by Remark \ref{r:all-is flat-and-lots-of-length}  
\begin{gather*}
\pi(U_{Q'^j}) \cap  I_1=\emptyset
\end{gather*}
where
\begin{gather*}
I_1:=\{x\in I_\eta: \sharp \pi^{-1}(x) = 1\}.
\end{gather*}
Hence, if
\begin{gather*}
U_{Q'^j} \cap (1-2^{-J+2}-4\epsilon_0)c_0Q' \neq \emptyset,
\end{gather*} 
then 
\begin{gather}\label{U_Q-j-gives-a-lot}
\pi(U_{Q'^j}) \cap  I_1=\emptyset.
\end{gather}
Now,
\begin{eqnarray*}
&&{\diam(U_{Q'}) \over
(1 + 2^{-J+1})(1-2^{-J+2}-4\epsilon_0)^{-1} (1-8\epsilon_0)^{-1}}\\
&\leq& {\diam (\eta_{Q'}) \over 
		(1 - 2^{-J+2}-4\epsilon_0)^{-1} }\\ 
&\leq& 
	\int\limits_{ I_{\eta}\cap (1-2^{-J+2}-4\epsilon_0)c_0Q'}  1\\
&\leq& 
	(\int\limits_{ (1-2^{-J+2}-4\epsilon_0)c_0Q' \cap I_1} 1+
	\int\limits_{ (1-2^{-J+2}-4\epsilon_0)c_0Q' \cap (I_\eta \smallsetminus I_1)} 1)\\
&\leq& {1\over 2}
	(\int\limits_{ (1-2^{-J+2}-4\epsilon_0)c_0Q' \cap I_1} 2+
	\int\limits_{ (1-2^{-J+2}-4\epsilon_0)c_0Q' \cap (I_\eta \smallsetminus I_1)} 2)\\
&\leq& {1\over 2} s'\,,
\end{eqnarray*}
where the last inequality follows from \eqref{e:2_in_R_Q} and \eqref{U_Q-j-gives-a-lot}.
Take $q<1$ such that 
\begin{gather*}
(1 + 2^{-J+1})(1-2^{-J+2}-4\epsilon_0)^{-1} (1-8\epsilon_0)^{-1}{1\over 2}\leq q
\end{gather*}
by enlarging $J$ and reducing $\epsilon_0$ if need be.

\noindent{\bf Step 2} and {\bf Step 3} are as in Lemma \ref{delta_1} replacing 
$U^\cfour$ with $U$.
\end{proof}

\begin{cor}
\begin{eqnarray*}
\sum\limits_{Q \in \Delta_{2.1}}\beta_{S_{Q}}(Q)\diam(Q)&\lesssim&
\length (\gamma)
\end{eqnarray*}
\end{cor}
%
%
%
\subsubsection{Putting it all together}
We now have
\begin{eqnarray*}
\sum\limits_{Q \in \G_2}\beta_{S_{Q}}(Q)\diam(Q) 
&\leq& 
\sum\limits_{Q \in \Delta_2}\beta_{S_{Q}}(Q)\diam(Q) + 
\sum\limits_{Q \in \Delta_1}\beta_{S_{Q}}(Q)\diam(Q)  \\
&\leq&
\sum\limits_{Q \in \Delta_{2.1}}\beta_{S_{Q}}(Q)\diam(Q) + 
\sum\limits_{Q \in \Delta_{2.2}}\beta_{S_{Q}}(Q)\diam(Q) + 
\sum\limits_{Q \in \Delta_1}\beta_{S_{Q}}(Q)\diam(Q)  \\
&\lesssim&
\length(\gamma)\\
&\lesssim&
\cH^1(\Gamma)
\end{eqnarray*}
We are done since we have shown  equation \eqref{sum_straight_eq}
\qed

This concludes the proof  of Theorem \ref{sum_beta_less_length} (and thus of Theorem \ref{int_beta_less_length} as well).
\begin{rem}
If one follows the computations one gets that the total constant here is dominated by 
$\sim A^{7\over 2}\log^2A$ which comes from $\Delta_{2.2}$.
\end{rem}

\section{Proof of Theorem \ref{construction_thm}}\label{constructions}
\subsection*{Farthest Insertion - A Local Version}\label{farthest_local_hilbert}

 
The following is a variation of what appears in \cite{J1}.  
Theorem \ref{construction_thm} can also be deduced as a special case from the work in \cite{Ha}.  We include a proof for completeness.

Let $K \subset H$  and $X^K=\cup X^K_n$ be given.
We  construct a connected set $\Gamma_0$ containing $K$ such that
\begin{gather}\label{100905}
\cH^1(\Gamma_0) \lesssim \diam(K) + \suml_{Q\in\hat{\G}^K} \beta_K^2(AQ)\diam(Q).
\end{gather}
By rearranging the constants ($A\to A^2$) we then get Theorem \ref{construction_thm}. 
We do this via  minor variations of the construction in \cite{J1}, which in turn is based on `Farthest Insertion' (see \cite{JM}) applied to the MST (Minimal Spanning Tree) problem. The only innovation here beyond \cite{J1}  is the treatment of the case $\beta\geq\epsilon_0$.  All other cases are treated in the same spirit.

If the right hand side of the inequality \eqref{100905} is infinite, then any connected set containing $K$ will suffice, and so we assume it is finite.
We may assume $K$ is closed without loss of generality.
We then get that   $K$ is compact. To see this, consider the contrary. 
We then have an infinite $\delta-net$,  $\{a_i\}$.  
Since $\diam(K)<\infty$,  we may assume WLOG that  there is no  infinite $2\delta-net$.
By perturbing, we may assume WLOG that $X^K_{\log(\delta)+3} \supset \{a_i\}$  and there is no infinite $3\delta-net$. 
Hence there is an $x\in X^K$ such that   $\ball(x,3\delta) \cap X^K_{\log(\delta)+3}$ is infinite.
If $A$ is large enough then this is a contradiction to $\sum_{Q \in \G} \beta(AQ)^2\diam(Q) < \infty$. 
Hence we will assume $K$ is compact.

We assume an order on $X^K$, such that  all points  in $X^K_n$ come before (are smaller than) points in $X^K_{n+1}\smallsetminus X^K_n$. We write $p_1< p_2<....$ where  $X=\{p_1,p_2,...\}$.
We define
\begin{gather*}
O(p_i):=\{p_j:j<i\}
\end{gather*}
and 
\begin{gather*}
d_i:=\dist(p_i,O(p_i)).
\end{gather*}


The construction below uses the following scheme.  
We inductively construct   a sequence of graphs $G_i$ with vertices in $X^K$. 
We freely confuse the graph $G_i$ with the set underlying the edges+vertices.
We will also have  a `virtual graph' $H_i$ which will simply be an addition to $G_i$ which will be used as an accounting tool. 
Start with the segment $[p_1,p_2]$.  Call it $G_2$.
Now, inductively, obtain $G_{i}$ from $G_{i-1}$ by connecting $p_i$. 
This may involve  connecting it to $G_{i-1}$ by modifying an edge that has both endpoints inside $\ball(p_i,A\cdot d_i)$ or by adding a new edge (the `cheaper' of the two options in a sense which will be clear later). Since changes are done only in $\ball(p_i,A\cdot d_i)$  
we refer to this as {\it a local version}.
In some cases we will  perform some preemptive constructions.  

We will see that the length is controlled by 
\begin{gather*}
\diam(K) +\suml_{i}\beta^2(p_i,A^2 2^{k(i)})A^2 2^{k(i)} 
\end{gather*}
where $p_i\in X_{k(i)}\smallsetminus X_{k(i)-1}$.
We will  get estimates for  $\cH^1(G_n \cup H_n)$.  
We set $H_2=[p_1,p_1+A(p_1-p_2)] \cup [p_2,p_2+A(p_2-p_1)]$.\\
We introduce some more notation:\\
At the induction stage we will add the point $x_0=p_{n_1}\in X^K_k\smallsetminus X^K_{k-1}$. 
We will call $Q=Q_{x_0}=B(x_0,A2^{k})$.  
If $x_0 \in G_{n_1-1}$ (underlying set) then we do nothing.  
Otherwise, for convenience of notation, denote the nearest point to it in 
$O(x_0)$ by $0$ (the origin). 
Denote by  $\R$ the line containing $x_0$ and passing through 
the origin, such that $x_0>0$. 
We will call $\pi_\R:H \to \R$ the orthogonal projection onto $\R$.
We will also write $\Re(z)=\pi_\R(z)$, borrowing notation from complex variable. 
Let\\
\indent$W=W_\R=\{z\in \half Q:  \frac{-\Re(z)}{\dist(z,\R)}\leq \frac{1}{\sqrt {3}} \}
\smallsetminus\{0\}$.\\ 
\indent$W^*=W^*_\R=\{z\in \half Q: \frac{\Re(z)}{\dist(z,\R)}\leq \frac{1}{\sqrt {3}} \}
\smallsetminus\{0\}$.\\     %
We assume that the following properties hold (we say that they hold {\it at $0$}):
\begin{verse}
(P1) If $x_0 \notin G_{n_1-1}$ and $\beta(Q) \leq \epsilon_0$, and $O(x_0) \cap W\neq \emptyset$, 
then
let $y_1 \in W$ minimize $\norm{z}$ on $O(x_0) \cap W$. 
We have $G_{n_1-1}\cup H_{n_1-1} \supset [0,y_1]$\\
(P2) If $x_0 \notin G_{n_1-1}$, $\beta(Q) \leq \epsilon_0$, 
and if $O(x_0) \cap W = \emptyset$  then 
there exists an interval $I \subset (G_{n_1-1}\cup H_{n_1-1}) \cap W  \st I=[0,z], z \in \partial (\half Q)$.\\
(P3) If $x_0 \notin G_{n_1-1}$, $\beta(Q) \leq \epsilon_0$, 
and if $O(x_0) \cap W^* = \emptyset$  then 
there exists an interval $I \subset (G_{n_1-1}\cup H_{n_1-1}) \cap W^*  \st I=[0,z], z \in \partial (\half Q)$.\\
\end{verse}
The way in which we will assure the assumptions   (P2) and (P3) will be by having segments starting from $0$ which have an angle between them which is large enough.

\begin{figure}[p]
\begin{center}
$\begin{array}{c}
\scalebox{0.3}{\includegraphics*[-2in,1in][14in,10in]{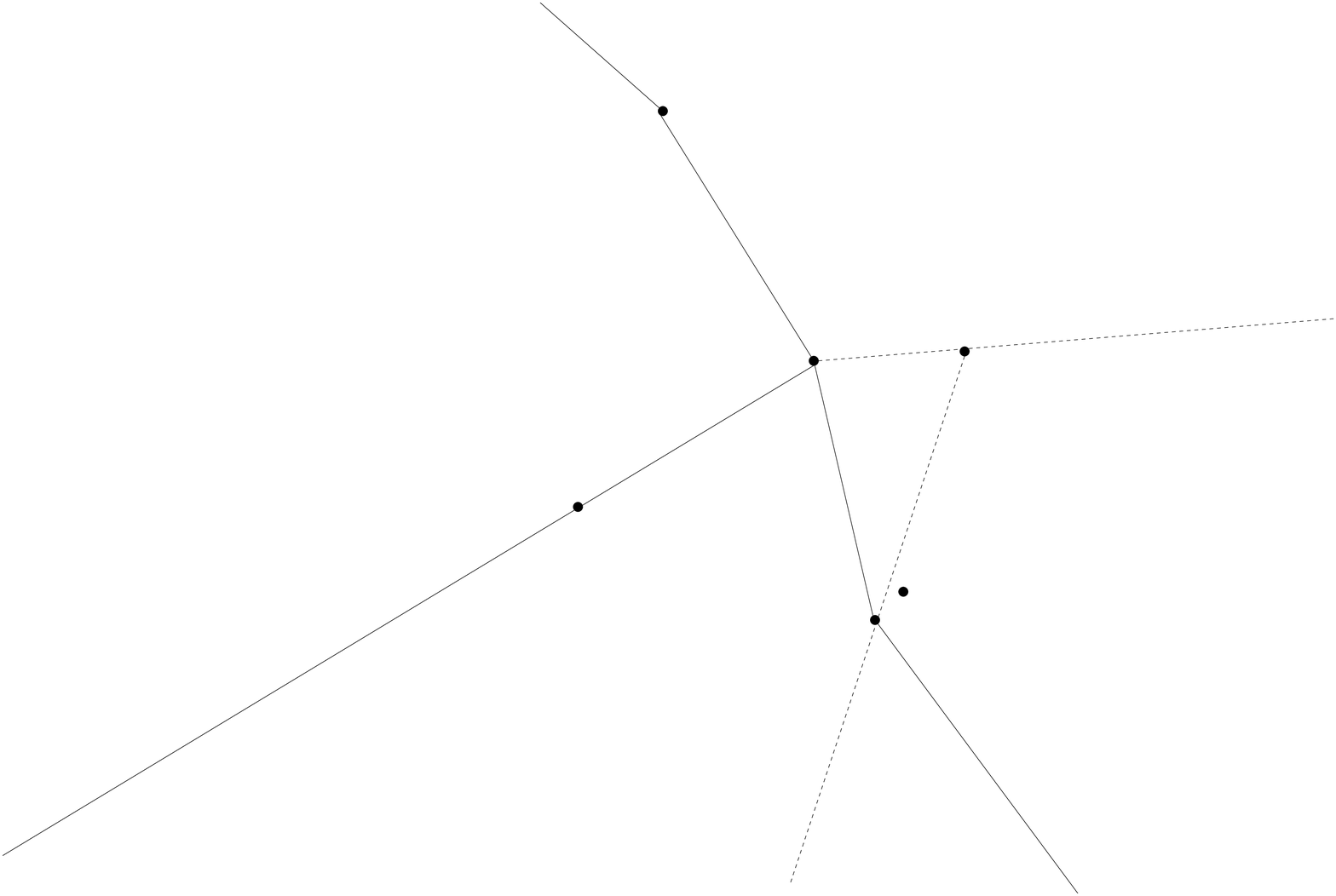}}
\put(-375,150){case 1:}
\put(-25,137){\small  $x_0$}
\put(-80,130){\small  $0$}\\
\scalebox{0.3}{\includegraphics*[-2in,0in][14in,4in]{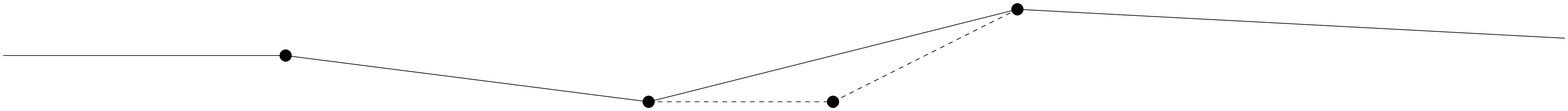}}
\put(-375,20){case 2:}
\put(-240,20){\small  $y_{-1}$}
\put(-60,15){\small  $y_1$}
\put(-105,-5){\small  $x_0$}
\put(-160,5){\small  $0$}\\
\scalebox{0.3}{\includegraphics*[-2in,-0in][14in,4in]{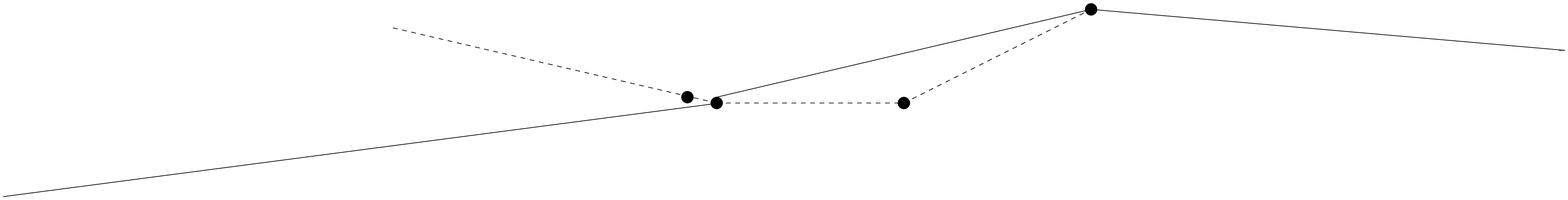}}
\put(-375,40){case 3:}
\put(-60,35){\small  $y_1$}
\put(-100,10){\small  $x_0$}
\put(-200,0){\small  $I$}
\put(-139,26){\small  $0$}\\
\scalebox{0.3}{\includegraphics*[-2in,0in][14in,4in]{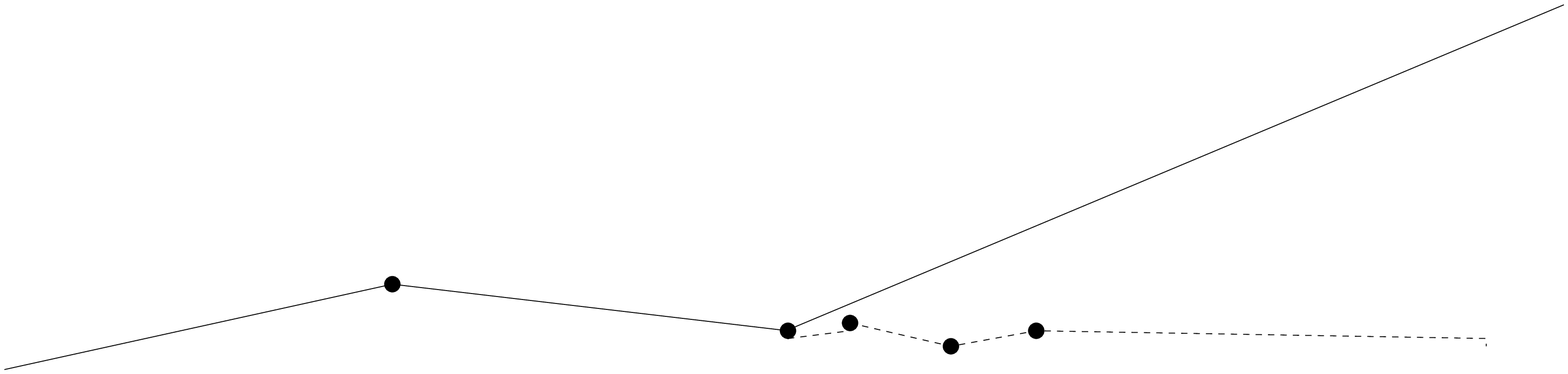}}
\put(-375,40){case 4:}
\put(-240,20){\small  $y_{-1}$}
\put(-126,10){\small  $x_0$}
\put(-80,50){\small  $I$}
\put(-170,10){\small  $0$}\\
\scalebox{0.3}{\includegraphics*[-2in,1.5in][14in,5in]{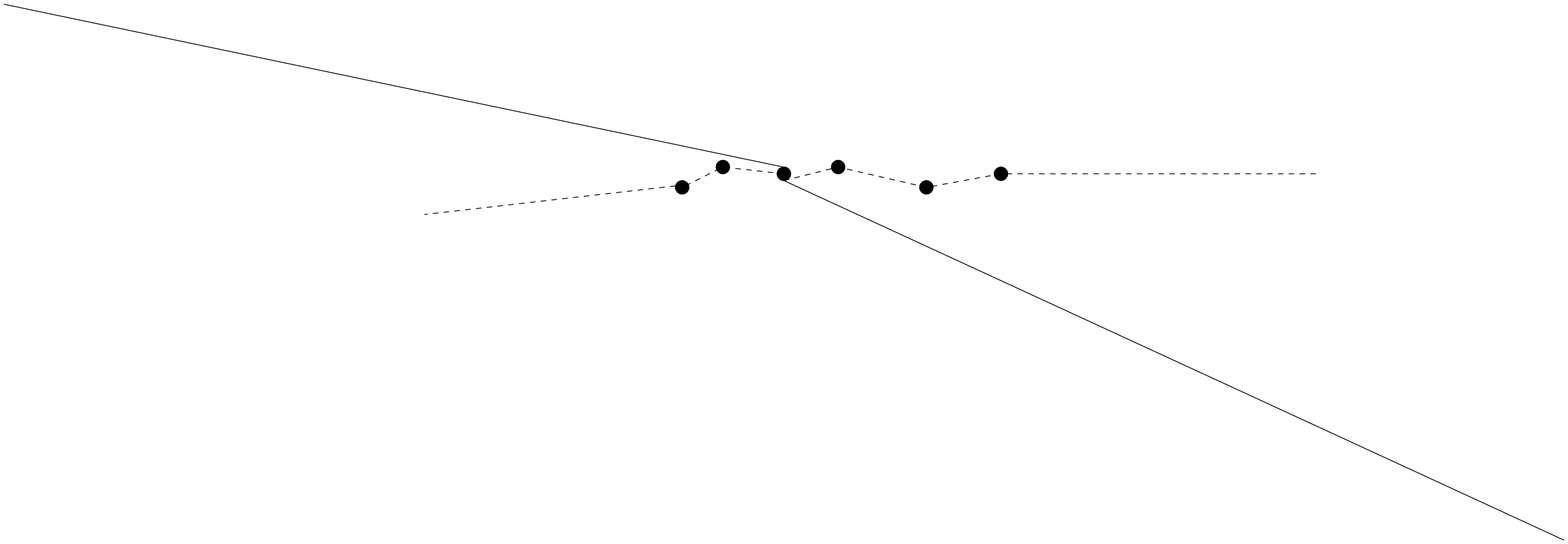}}
\put(-375,30){case 5:}
\put(-105,45){\small  $x_0$}
\put(-220,60){\small  $I$}
\put(-90,0){\small  $I$}
\put(-150,45){\small  $0$}\\
\end{array}$
\end{center}
\caption{Dotted lines are additions to $G_{n_1-1}\cup H_{n_1-1}$ giving $G_{n_1}\cup H_{n_1}$}
\label{construction-cases}
\end{figure}

We now go over the possible cases.  See Figure \ref{construction-cases} for examples.

\textbf{Case 1:} 
	$\beta(Q) \geq \epsilon_0$.\\
Add  segments in order to most efficiently connect $x_0$ to $G_{n_1-1}$ by connecting to  vertices in 
$O(x_0) \cap Q$.
This can involve adding a segment connecting $x_0$ to a vertex in $O(x_0) \cap Q$, or 
modifying a segment in $G_{n_1-1}$.   
Further, we make sure (P2) and (P3) are  preserved at $x_0$ by adding  to $H_{n_1}$ segments 
(if necessary) at $x_0$.
The costs of the addition to $G_{n_1-1}\cup H_{n_1-1}$ is a total of $\lesssim\beta(Q)^2\diam(Q)$.
To make sure (P1) is preserved we  need only consider  points in $X^K_{k+\log(A)}$.  
For each of those we  make sure (P1) is preserved by adding the required segment.  
We may do so for \textbf{ALL} points of $X^K_k$ at a total cost of 
$\lesssim\suml_{Q=B(x,A2^{-k})\atop x\in X^K_{k+\log(A)}}\beta(AQ)^2\diam(Q)$.
This is the key difference with \cite{J1}, where Jones used local compactness of $\R^d$ to give bounds on the length added, and was hence able to be more `wasteful' in adding segments.

\textbf{Case 2:} 
	$\beta(Q) < \epsilon_0; 
	O(x_0) \cap W \neq \emptyset; 
	O(x_0) \cap W^* \neq \emptyset$.\\
Let $y_1 \in W$ minimize $\norm{z}$ on $O(x_0) \cap W$.
By P1 we have $[0,y_1] \subset G_{n_1-1}\cup H_{n_1-1}$.  
If $[0,y_1] \subset G_{n_1-1}$ we replace it with $[0,x_0],[x_0,y_1]$ in  $G_{n_1}$.\\ 
If $[0,y_1] \subset H_{n_1-1}$ then we replace it with $[0,x_0],[x_0,y_1]$, 
placing the shorter of the above intervals in $G_{n_1}$ and the longer one in $H_{n_1}$.\\ 
We  get that (P1), (P2) and (P3)  are maintained.
By the Pythagorean theorem (and a first order approximation) we get that the length added 
to $G_{n_1-1}\cup H_{n_1-1}$ is bounded by $C\beta(Q)^2\diam(Q)$.  

\textbf{Case 3:} 
	$\beta(Q) < \epsilon_0$; 
	$O(x_0) \cap W \neq \emptyset$;
	$O(x_0) \cap W^* = \emptyset$.\\
This corresponds to Cases 3 and 4 in \cite{J1}.
Let $y_1 \in W$ minimize $\norm{z}$ on $O(x_0) \cap W$.
If $[0,y_1] \subset G_{n_1-1}$ then replace it with $[0,x_0],[x_0,y_1]$ in  $G_{n_1}$.\\ 
If $[0,y_1] \subset H_{n_1-1}$ then replace it with $[0,x_0],[x_0,y_1]$, placing the shorter of the above intervals in $G_{n_1}$ and the longer one in $H_{n_1}$.\\
(P1), (P2) and (P3) are maintained at $x_0$ since $\beta<\epsilon_0$.\\ 
We now make sure they are maintained at $0$ and $y_1$.
Denote by $z_{1}$ the point maximizing $\norm{z}$ in $(\{0\}\cup W^*)\cap X^K_{k+\log(A)}$.  
Since $W^*\cap O(x_0)=\emptyset$ we have $\norm{z_1}\leq 2^{-k+1}$.  
Let $z_1,...,z_N=y_1$ be the points in $\ball(0,\norm{y_1})\cap X^K_{k+\log(A)}$,  ordered by increasing `Real' value $\Re(\cdot)$.
\indent Add on to $G_n$ the segments $[z_1,z_2],...,[z_{N-1},z_N]$.\\
\indent 
	If $\norm{z_{1}} \neq 0$  add on to $H_n$ the segment 
		$[z_1,2^k {z_1\over\norm{z_1}}]$.  
	Otherwise, add on to $H_n$ the segment $[0,-2^{-k}]$.\\
If $y_1$ also maximizes $\norm{z}$ on $O(x_0) \cap W$ then do a similar (symmetric) construction near $y_1$.  Otherwise do nothing there.\\
The point is that the segments we have added will give us (P1), (P2), and (P3) $\log(A)$ scales into the future, and until then we have done all constructions needed.

We now need to account for the length we have added. 
(P2) assures us a line segment in the form of $[0,z]\subset G_{n_1-1} \cup H_{n_1-1}$ where 
$\norm{z}=A2^{-k-1},  z\in Q$.
Set $I_Q=[A2^{-k-4},A2^{-k-3}]$.
$I_Q$ will not be altered (or moved) at any future stage (as it is far from $K$ and deep inside $W^*$).  
Hence $I_Q \subset \bigcap\limits_{k=0}^\infty (G_{n_1-1+k}\cup H_{n_1-1+k})$.  
The length that we added is bounded by 
\begin{gather*}
C\beta^2(Q)\diam(Q) + 2\cdot2^{-k} +3\cdot2^{-k+1}({1\over A} + 2\epsilon_0 A)A 
\leq {100 \over A} \cH^1(I_Q)
\end{gather*}
by reducing $\epsilon_0$.
We also note that $x \in I_Q$ will be used as such at most once.\\

%
%
\textbf{Case 4:} 
	$\beta(Q) < \epsilon_0$; 
	$O(x_0) \cap W = \emptyset$;
	$O(x_0) \cap W^* \neq \emptyset$\\
This corresponds to Case 5 in \cite{J1}, and is similar in accounting to Case 3 above.
Let $y_{-1} \in W$ minimize $\norm{z}$ on $O(x_0) \cap W^*$.
We have $\norm{x_0}\leq \norm{y_{-1}}$.  
Let $z_1,...,z_N$ be the points in $(\{0\}\cup W)\cap X^K_{k+\log(A)}$,  ordered by increasing `Real' value $\Re(\cdot)$.
Add on to $G_n$ the segments\\
$[z_1,z_2],...,[z_{N-1},z_N]$. 
Add on to $H_n$ the segment 
		$[z_N,2 z_N]$. 

\textbf{Case 5:} 
	$\beta(Q) < \epsilon_0$; 
	$O(x_0) \cap W = \emptyset$;
	$O(x_0) \cap W^* = \emptyset$\\
This corresponds to Case 6 in \cite{J1}, and is similar in accounting to Case 3 and 4 above.
We have $\norm{x_0}\leq 2^{-k+1}$.  
Let $z_1,...,z_N$ be the points in $\ball(0,2^{-k+1})\cap X^K_{k+\log(A)}$,  ordered by increasing `Real' value $\Re(\cdot)$.
Add on to $G_n$ the segments $[z_N,z_{N-1}],...,[z_2,z_1]$.
Add on to $H_n$ the segment 
		$[z_N,2 z_N]$. 
If $\norm{z_{1}} \neq 0$  add on to $H_n$ the segment 
		$[z_1,2^k {z_1\over\norm{z_1}}]$.
	Otherwise, add on to $H_n$ the segment $[0,-2^{-k}]$.\\


This concludes all the cases.

Inductively we get:
\begin{eqnarray*}
&&\cH^1(G_n\cup H_n)-\cH^1(G_2\cup H_2)\\
 &\leq&
\suml_{x\in O(p_{n})} [C\beta(Q_x)^2\diam(Q_x)  + {100 \over A} \cH^1(I_{Q_x})]  
+ C\log(A)\suml_{Q\in \hat{\G}}\beta(AQ)^2\diam(Q)\\
&\leq&
C\sum_{Q \in \hat{\G}} \beta(Q)^2\diam(Q) 
+ C\suml_{Q\in \hat{\G}}\beta(AQ)^2\diam(Q) 
+      \suml_{Q=Q_x:x\in O(p_{n})}{100 \over A} \cH^1(I_Q)\\
&\leq&
C\sum_{Q \in \hat{\G}} \beta(Q)^2\diam(Q) + 
C\suml_{Q\in\hat{\G}}\beta(AQ)^2\diam(Q)+ 
{100 \over A} \cH^1(G_{n}\cup H_n)\\
&\leq&
C\suml_{Q\in\hat{\G}}\beta(AQ)^2\diam(Q)+ 
{100 \over A} \cH^1(G_{n}\cup H_n).\\
\end{eqnarray*}
Hence 
\begin{gather*}
\cH^1(G_{n})\leq\cH^1(G_{n}\cup H_n) \lesssim \sum_{Q \in \hat{\G}} \beta(AQ)^2\diam(Q) + \diam(K)
\end{gather*} 
by choosing $A$ large enough.  Denote such a choice by $A_0$.\\

Set 
\begin{gather*}
E:=\{x\in H:x=tx_1 + (1-t)x_2, x_i \in K', -A\leq t \leq A\}.
\end{gather*}
We get $E$ is compact by Lemma \ref{cpt}.
We also have $G_n \subset E$.
We use Lemma \ref{parametrization}
to get $\gamma_n:[0,1] \to G_n$ Lipschitz with uniformly bounded Lipschitz norm.
Since $E$ is compact we can use Arzela-Ascoli to obtain a limit $\gamma$.  
We have 
\begin{gather*}
K^{closure} \subset Image(\gamma)=:\Gamma_0
\end{gather*}
and we get the desired estimate on $\cH^1(\Gamma_0)$.\\

\qed

Note that the exact same computation works for $G_n\cup H_n$.

\begin{rem}\label{FI}
One may take $X^K_1:=\{p_1,p_2\}$ with $\dist(p_1,p_2)=\diam(K)$ and inductively or $i>2$, $p_i$ maximizing $d_i$, and define $X^K_n:=\{p_i:d_i<2^{-n}\}$ and  take the order given by the induction above in each $X^K_n$.  Then aside from the preemptive constructions appearing in Cases 3-5, our algorithm gives a {\it local} version of the {\bf Farthest Insertion} algorithm. Note however, that these preemptive constructions are exactly the constructions that would have appeared  under a {\it local} version of the Farthest Insertion algorithm, by which we mean connecting $p_n$ in the most efficient way to vertices in $O(p_n)\cap Q_{p_n}$.
See \cite{JM} for the standard version.
\end{rem}

%
%
%
%
%
%
%
%
%
%
%
%
%
%
%
%
%
%
%
%
%
%
%
%

%

\section{Appendix}
%
%
\subsection{Proofs of Point-Set Topology Lemmas}\label{appendix_point-set}
\smallskip
\noindent
{\bf Lemma \ref{closure-length}}
{\it
Assume $\Gamma$ is connected.  Then $\cH^1(\Gamma)=\cH^1(\Gamma^{closure})$.
}
\begin{proof}
First note that since $\cH^1(\Gamma)\leq\cH^1(\Gamma^{closure})$ we need only concern ourselves with the case $\cH^1(\Gamma)<\infty$.

Let $\epsilon>0$ and $\delta>0$ be given.  By the definition of Hausdorff measure, there exists sets $\{E_i\}_0^\infty$ such that 
$\sum\diam(E_i)\leq\cH^1(\Gamma)(1+\epsilon)$ and 
$\diam(E_i)\leq \delta$.  
Without loss of generality we may assume  the sets $E_i$ are convex and hence connected.
We may also assume without loss of generality that they are open 
(by first taking sets that have diameters summing up to $\cH^1(\Gamma)(1+\half\epsilon)$ and then taking small neighborhoods of them).
One may construct from them families $\{G_j\}_0^\infty$ and $\{B_j\}_0^\infty$ such that 

(i)   $\delta\leq\diam(G_j)<2\delta$ and $\diam(B_j)<\delta$

(ii)  both $G_j$ and $B_j$ are  connected

(iii) both $G_j$ and $B_j$ are  unions of some sets $E_i$ 

(iv) $\bigcup G_j \cup \bigcup B_j=\cup E_i$

(v) $B_j\cap B_k=\emptyset$  for all $j\neq k$

(vi)  if $B_j\cap G_k\neq \emptyset$ then $\diam(B_j\cup G_k)\geq 2\delta$.

This can be done by inductively going over the sets $E_i$ and joining them whenever possible.
Since $\Gamma$ is connected, every set $B_j$ intersects some set $G_i$ (possibly more then one).  Denote a choice of such a set by $G_{i(j)}$.  
Set $F_i=G_i\cup \bigcup\limits_{\{j: i(j)=i\}}B_j$.
We get by the triangle inequality

(i) $\delta \leq \diam(F_i)< 4\delta$

(ii) $\sum\diam F_i \leq \sum\diam(E_i) \leq \cH^1(\Gamma)(1+\epsilon)$
 
We conclude that there is a finite number of sets $F_i$ and we may consider a $\half\epsilon\delta$ neighborhood of them, $F'_i$.  
We have 
\begin{gather*}
\cup F'_i\supset \Gamma^{closure}
\end{gather*}
We also have
\begin{gather*}
\diam (F'_i) 
\leq diam(F_i)  + \epsilon\delta 
\leq diam(F_i) (1 + \epsilon) 
\end{gather*}
and so 
\begin{gather*}
\sum\diam(F'_i)\leq \cH^1(\Gamma)(1+\epsilon)^2
\end{gather*}
We conclude that $\cH^1(\Gamma^{closure})\leq \cH^1(\Gamma)$ and so we have equality as desired.
\end{proof}

\smallskip
\noindent
{\bf Lemma \ref{finite_length_then_cpt}}
{\it
Assume $\Gamma\subset H$ is a closed connected set with $\cH^1(\Gamma)< \infty$.
Then $\Gamma$ is compact.
}
\begin{proof}
Assume $\Gamma$ is not compact.  Hence for arbitrarily small $\delta$ we can obtain an infinite $\delta-net$ for $\Gamma$: $\{a^\delta_n\}$.  
By connectedness we have $\cH^1(\Gamma \cap \ball(a^\delta_n,\third\delta)) \geq\third\delta$.
We also have $\ball(a^\delta_n,\third\delta)$ are disjoint, which contradicts $\cH^1(\Gamma)< \infty$.
\end{proof}

\smallskip
\noindent
{\bf Lemma \ref{cpt}}
{\it
Let $C_1,C_2 >0$ be given.
Given a compact connected set $\Gamma \subset H$ the set
$E:=\{x\in H:x=tx_1 + (1-t)x_2, x_i \in \Gamma, -C_1\leq t \leq C_2\}$
is compact.
}
\begin{proof}
Suppose $\{x^i\} \subset E$ is a sequence.  
We can write $x^i=t^ix^i_1 + (1-t^i)x^i_2$ as in the definition of $E$.
By the compactness of $\Gamma$ we have 
$i_k \st \\
x^{i_k}_1 \to x_1$. By compactness of $\Gamma$ again, $x^{i_{k_j}}_2 \to x_2$. 
By  compactness of  $[-C_1,C_2]$ we have $t^{i_{k_{j_l}}} \to t$. 
$x_1,x_2 \in \Gamma, t \in [-C_1,C_2]$.  Hence $x^{i_{k_{j_l}}}\to tx_1+ (1-t)x_2 \in E$.
\end{proof}

One can find an $\R^d$ version of the following lemma in \cite{DS} with the 
Lipschitz norm depending on $d$.
The following proof is a modification of the proof given there, which gives a result 
independent  of $d$.

\smallskip
\noindent
{\bf Lemma \ref{parametrization}}
{\it
Let $\Gamma \subset H$ be a compact connected  set of finite length.  
Then we have a Lipschitz function $\gamma:[0,1] \to H \st Image (\gamma)=\Gamma$ 
and $\norm{\gamma}_{Lip} \leq 32\cH^1 (\Gamma)$
}
\begin{proof}
We use a  well known result from graph theory:\\
If $G$ is a connected graph with finitely many edges, then there is a path that traverses each 
edge of $G$ exactly twice (once in each direction).
This result is easily seen by induction on the number of edges.

For $n\geq0$, let $X_n=X_n^\Gamma$ (i.e. take $X_n \subset \Gamma$ a $2^{-n} - net$  such that 
$X_n \subset X_{n+1}$).\\
We want to get a connected set $E_n$.
We do this by adding line segments inductively.
Set $E_n^0=X_n$.  
We get  $E_n^{i+1}$ form $E_n^i$ by adding a line segment between points 
$x_1,x_2 \in X_n \st \dist(x_1,x_2)<2^{-n+3}$ and they are not yet in the same 
connected component of $E_n^i$.  
If there are no two such points we stop and call the resulting set $E_n$.
Let $G_n$ be the obvious abstract graph associated to $E_n$. 
If $G_n$ is not connected then ${\rm Vertex} (G_n)=A \cup B$ with $\dist(A,B) \geq 2^{-n+2}$ and  
$A$ separated from $B$.
By the construction of $E_n$ and $X_n$ we have that $\dist(\N_{2^{-n}}(A),\N_{2^{-n}}(B)) \geq 2^{-n}$ and $\Gamma \subset \N_{2^{-n}}(A)\cup \N_{2^{-n}}(B)$.  This is a contradiction to  $\Gamma$ being connected.  
Hence $G_n$ is connected.\\
Note that $\cH^1(E_n) \leq \sharp(X_n)2^{-n+3} \leq 16 \cH^1(\Gamma)$, where the final inequality follows from the fact that the balls $\{B(x,2^{-n-1}):x \in X_n\}$ are disjoint.\\
We can thus parameterize $E_n$ by a Lipschitz curve of $\gamma_n:[0,1] \to H$. 
The image of this parameterization is in $E$ as defined in the previous lemma.  
By Arzela-Ascoli we have a subsequence converging to $\gamma$.
We have that $Image(\gamma)=\Gamma$ by say 
\begin{gather*}
\supl_{x\in E_n} \dist(x,\Gamma)+ \supl_{y\in \Gamma} \dist(E_n,y) \leq 3\cdot 2^{-n} +   2^{-n}
= 4\cdot 2^{-n}
\end{gather*}
 and a triangle inequality. 
\end{proof}

Hence, we also have

\smallskip
\noindent
{\bf Corollary \ref{cor-parametrization}}
{\it
Let $\Gamma \subset H$ be a compact connected  set of finite length.  
Then we have a Lipschitz function $\gamma:\T \to H \st Image (\gamma)=\Gamma$ 
and $\norm{\gamma}_{Lip} \leq 32\cH^1 (\Gamma)$
}

%
%
%
%
%
%
%
%
%
%
%
%
%
%
%
%
%
%
%
%
%
%
%
%
%
%
\subsection{Existence of the MST}\label{appendix_MST}
\begin{lemma}\label{MST-existence}
Assume $K\subset \Gamma_0$,
where $\Gamma_0$ is the image of an arc-length  parameterization $\gamma_0$ 
and $\cH^1(\Gamma_0)<\infty$.
Then there exists a curve with image $\Gamma_{MST}\supset K$ such that 
$\cH^1(\Gamma)\geq \cH^1(\Gamma_{MST}),\forall \Gamma \supset K$ 
\end{lemma}
\begin{proof}
We have $E:=\{x\in H:x=tx_1 + (1-t)x_2, x_i \in \Gamma_0, 0\leq t \leq 1\}$ is compact by Lemma \ref{cpt}.
Consider also $F=K^{\rm closure}$ which is also compact as a closed subset of a compact set.

Assume that $\gamma_n$ are such that $\cH^1(\Gamma_n)$ decreases to 
$L=\inf_{\Gamma\supset K}\cH^1(\Gamma)$, where $\Gamma_n$ is the image of $\gamma_n$.
WLOG we may assume that for all $n$, $\gamma_n$ is defined on $\T$ and is Lipschitz with constant $C\leq32(L+1)$.
We will define $\gamma^*_n$ as follows:
\begin{quote}
Assume $[a_i,b_i]\subset\T$ is a maximal interval with an interior whose image under $\gamma_n$ is disjoint from $F$. 
Set $\gamma^*_n 1_{[a_i,b_i]}=\pi_{<\gamma_n(a_i),\gamma_n(b_i)>}\gamma_n$, 
where \newline
$\pi_{<\gamma_n(a_i),\gamma_n(b_i)>}$ 
is projection onto the line going 
through 
$\gamma_n(a_i)$ and $\gamma_n(b_i)$.If $\gamma_n(a_i)=\gamma_n(b_i)$ define $\pi_{<\gamma_n(a_i),\gamma_n(b_i)>}=\gamma_n(a_i)$.
\end{quote}
We get that $\gamma^*_n$ is continuous and has a derivative a.e. bounded by $C$ (as the set of exit/entry points from $E$ is countable). We also have that the image of every $\gamma^*_n$ is in $E$.
Now use Arzela-Ascoli to get a limit path $\gamma$ with image $\Gamma_{MST}$ (whose name is yet to be justified).  WLOG assume that the original sequence $\gamma_n$ converges to $\gamma$.
Using Golab's Theorem (for subsets of $\R^d$) we have
\begin{gather*}
L\leq\cH^1(\Gamma_{MST})
=\lim_d(\cH^1(\pi_d\Gamma_{MST}))
\leq\lim_d(\liminf\cH^1(\pi_d\Gamma_n))
\leq \lim_d (\liminf\cH^1(\Gamma_n))
=L
\end{gather*}
and so  $L=\cH^1(\Gamma_{MST})$ ($\pi_d$ is projection onto the first $d$ coordinates). (The first equality above follows for instance by an idea similar to the one of the proof of Lemma \ref{closure-length}.)
\end{proof}

%
%
%
%
%
%
%
%
%
%
%
%
%
%
%
%
%
%
%
%
%
%
%
%
%
%
\subsection{Table of Notation/Symbols}\label{appendix_notation-table}
Below is a table of notation/symbols.  For each element, we list the first time it is defined or mentioned (or repeat the definition if it is short).

\medskip

\begin{tabular}{l@{\extracolsep{.6in}}l}
{\bf Symbol}& {\bf Location}\\
\\
$A$ &  see  \eqref{G-K-def}.\\
$\bt$ &  see \eqref{bt-definition}.\\
$c_0$ &  following Lemma \ref{building_U}.\\
$C_U$ &  see Remark \ref{09042005} and preceding definition.\\
$\Delta_1,\ \Delta_2,\ \Delta_{2.1},\  \Delta_{2.2}$&  following Proposition \ref{sum_straight_prop}.\\
$\epsilon_1$& definition of $\G_i^j$ (see below).\\
$\epsilon_2$& following \eqref{bt-definition}.\\
$\F$& Section \ref{curvy_arcs}.  Either a {\it Filtration} (see Lemma \ref{L2_lemma} or \\
                 & a collection on-route to a {\it Filtration}.\\
$\gamma$ & following Corollary \ref{cor-parametrization}.\\
$\gamma_Q$ & following \eqref{S_Q-definition}.\\
$\Gamma$& Theorem \ref{sum_beta_less_length}.\\
$\G$& see \eqref{g-def}.\\
$\G_i^j$& following \eqref{S_Q-definition}.\\
$\hat{\G},\  \hat{\G^K}$& see \eqref{G-K-def}.\\
$J$& appears in several places as a factor in the number\\
& of scales {\bf J}umped.\\
$\ell$& arclength measure or its pushforward.\\
$\Lambda(Q) $& see \eqref{lambda-definition}.\\
$M$& see subsection \ref{delta-1-subsub}.\\
$\N_\epsilon(E)$& an $\epsilon$ neighborhood of $E$.\\
$Q$& a multiresolution element (e.g. $Q\in \G$).\\
$S_Q$& see \eqref{S_Q-definition}.\\
$U_Q, U_Q^x, U_Q^{xx}$& following Lemma \ref{building_U}.\\
$w_Q$& a weight/density.  Defined (differently) in the proofs of\\
& Lemma \ref{delta_1} and Lemma \ref{delta_2.1}.
\end{tabular}

\bibliographystyle{alpha}
\bibliography{../../../bibliography/bib-file}   

\begin{thebibliography}{Scharb}

\bibitem[Aro03]{Ar}
Sanjeev Arora.
\newblock Approximation schemes for {NP}-hard geometric optimization problems:
  a survey.
\newblock {\em Math. Program.}, 97(1-2, Ser. B):43--69, 2003.
\newblock ISMP, 2003 (Copenhagen).

\bibitem[BJ90]{BJ}
Christopher~J. Bishop and Peter~W. Jones.
\newblock Harmonic measure and arclength.
\newblock {\em Ann. of Math. (2)}, 132(3):511--547, 1990.

\bibitem[Chr90]{Ch}
Michael Christ.
\newblock A {$T(b)$} theorem with remarks on analytic capacity and the {C}auchy
  integral.
\newblock {\em Colloq. Math.}, 60/61(2):601--628, 1990.

\bibitem[Dav91]{Da}
Guy David.
\newblock {\em Wavelets and singular integrals on curves and surfaces}, volume
  1465 of {\em Lecture Notes in Mathematics}.
\newblock Springer-Verlag, Berlin, 1991.

\bibitem[DS93]{DS}
Guy David and Stephen Semmes.
\newblock {\em Analysis of and on uniformly rectifiable sets}, volume~38 of
  {\em Mathematical Surveys and Monographs}.
\newblock American Mathematical Society, Providence, RI, 1993.

\bibitem[FFPar]{FFP}
Fausto Ferrari, Bruno Franchi, and Herv{\'e} Pajot.
\newblock The geometric traveling salesman problem in the {H}eisenberg group.
\newblock {\em Rev. Mat. Iberoamericana}, To appear.

\bibitem[Hah05]{Ha}
Immo Hahlomaa.
\newblock Menger curvature and {L}ipschitz parametrizations in metric spaces.
\newblock {\em Fund. Math.}, 185(2):143--169, 2005.

\bibitem[Hahar]{Ha-2}
Immo Hahlomaa.
\newblock Curvature and {L}ipschitz parametrizations in 1-regular metric
  spaces.
\newblock {\em Annales Academiae Scientiarum Fennicae}, To appear.

\bibitem[JM02]{JM}
David~S. Johnson and Lyle~A. McGeoch.
\newblock Experimental analysis of heuristics for the {STSP}.
\newblock In {\em The traveling salesman problem and its variations}, volume~12
  of {\em Comb. Optim.}, pages 369--443. Kluwer Acad. Publ., Dordrecht, 2002.

\bibitem[Jon88]{J2}
Peter~W. Jones.
\newblock Lipschitz and bi-{L}ipschitz functions.
\newblock {\em Rev. Mat. Iberoamericana}, 4(1):115--121, 1988.

\bibitem[Jon90]{J1}
Peter~W. Jones.
\newblock Rectifiable sets and the traveling salesman problem.
\newblock {\em Invent. Math.}, 102(1):1--15, 1990.

\bibitem[KK92]{KK}
Claire Kenyon and Richard Kenyon.
\newblock How to take short cuts.
\newblock {\em Discrete Comput. Geom.}, 8(3):251--264, 1992.
\newblock ACM Symposium on Computational Geometry (North Conway, NH, 1991).

\bibitem[Ler03]{Le}
Gilad Lerman.
\newblock Quantifying curvelike structures of measures by using {$L\sb 2$}
  {J}ones quantities.
\newblock {\em Comm. Pure Appl. Math.}, 56(9):1294--1365, 2003.

\bibitem[Mat95]{Ma}
Pertti Mattila.
\newblock {\em Geometry of sets and measures in {E}uclidean spaces. {F}ractals
  and rectifiability}, volume~44 of {\em Cambridge Studies in Advanced
  Mathematics}.
\newblock Cambridge University Press, Cambridge, 1995.

\bibitem[Oki92]{Ok}
Kate Okikiolu.
\newblock Characterization of subsets of rectifiable curves in {${\bf R}\sp
  n$}.
\newblock {\em J. London Math. Soc. (2)}, 46(2):336--348, 1992.

\bibitem[Paj02]{Pa1}
Herv{\'e} Pajot.
\newblock {\em Analytic capacity, rectifiability, {M}enger curvature and the
  {C}auchy integral}, volume 1799 of {\em Lecture Notes in Mathematics}.
\newblock Springer-Verlag, Berlin, 2002.

\bibitem[Sch05]{Schul}
Raanan Schul.
\newblock {\em Subset of Rectifable curves in {H}ilbert Space and the Analyst's
  TSP}.
\newblock PhD thesis, Yale University, 2005.

\bibitem[Schara]{RS-metric}
Raanan Schul.
\newblock Ahlfors-regular curves in metric spaces.
\newblock {\em Annales Academiae Scientiarum Fennicae}, To appear.

\bibitem[Scharb]{my-TSP-survey}
Raanan Schul.
\newblock Analyst's traveling salesman theorems. {A} survey.
\newblock {\em Proceedings of the Ahlfors-Bers Colloquium}, To appear.

\bibitem[SS83]{StSt}
E.~M. Stein and J.-O. Str{\"o}mberg.
\newblock Behavior of maximal functions in {${\bf R}\sp{n}$} for large {$n$}.
\newblock {\em Ark. Mat.}, 21(2):259--269, 1983.

\bibitem[Ste83]{St}
E.~M. Stein.
\newblock Some results in harmonic analysis in {${\bf R}\sp{n}$}, for
  {$n\rightarrow \infty $}.
\newblock {\em Bull. Amer. Math. Soc. (N.S.)}, 9(1):71--73, 1983.

\end{thebibliography}

\end{document}